\documentclass [10pt,a4paper]{article}
\usepackage{amssymb}
\usepackage{enumerate,verbatim}
\usepackage{amsmath}
\usepackage{amsfonts,mathrsfs}
\usepackage{amsthm,wasysym}
\usepackage{epsfig,lpic,tikz}
\usepackage[applemac]{inputenc}
\usepackage[english]{babel}

\usepackage[width=0.9\textwidth]{caption}

\usepackage{pdfsync}
\usepackage{xcolor,hyperref}
\usepackage[T1]{fontenc}
\usepackage{lpic}
\usepackage{mathtools}

\hypersetup{
    colorlinks=true,
    linkcolor=blue,
    filecolor=magenta,      
    urlcolor=blue,
    pdftitle={Overleaf Example},
    pdfpagemode=FullScreen,
    }

\definecolor{green}{rgb}{0,0.5,0}

\usepackage{marvosym}

\newcounter{teoremaganso}
\newcounter{appendix}

\newcounter{coryganso}

\renewenvironment{abstract}{\small\quotation\noindent
 {\bfseries \abstractname .}}{\endquotation \par}

\newenvironment{trivlistparm}[2]{\trivlist
\item[{#1 #2}]}{\endtrivlist}
\newenvironment{proclama}[1]{\trivlistparm{\bfseries}{#1}\itshape}{\endtrivlistparm}
\newenvironment{prooftext}[1]{\trivlistparm{\bfseries}{#1}}{\Qed\endtrivlistparm}
\newenvironment{prova}{\trivlistparm{\bfseries}{Proof.}}{\Qed\endtrivlistparm}

\catcode`\@=11

\def\resetthefootnote{\renewcommand{\thefootnote}{\@arabic\c@footnote} }
\def\@principiremex#1{\trivlist
 \item[\hskip \labelsep{\bfseries #1\ \thetheo.}]\ignorespaces}
\def\opar@principiremex#1[#2]{\trivlist
 \item[\hskip \labelsep{\bfseries #1\ \thetheo\ (#2).}]\ignorespaces}

\newcommand{\newTHEOremrom}[2]{\newenvironment{#1}{\refstepcounter{theo}\@ifnextchar[{\opar@principiremex{#2}}
{\@principiremex{#2}}}{\qedB\endtrivlist}} \catcode`\@=12
\DeclareMathSymbol{\square}{\mathord}{AMSa}{"03}
\newcommand{\qedB}{\nopagebreak\hspace*{\fill}$\square$\par}
\newcommand{\Qed}{\nopagebreak\hspace*{\fill}{\vrule width6pt height6pt depth0pt}\par}

\newtheorem {theo} {Theorem} [section]
\newtheorem {prop} [theo] {Proposition}

\newtheorem {lem} [theo] {Lemma}
\newtheorem {bigtheo} [teoremaganso] {Theorem}

\newTHEOremrom {defi} {Definition}
\newTHEOremrom {obs} {Remark}
\newTHEOremrom {ex} {Example}

\newcommand{\refc}[1]{\mbox{$(\ref{#1})$}}
\newcommand{\secc}[1]{Section~\ref{#1}}

\newcommand{\teoc}[1]{Theorem~\ref{#1}}
\newcommand{\propc}[1]{Proposition~\ref{#1}}

\newcommand{\lemc}[1]{Lemma~\ref{#1}}
\newcommand{\defic}[1]{Definition~\ref{#1}}
\newcommand{\obsc}[1]{Remark~\ref{#1}}

\newcommand{\figc}[1]{Figure~\ref{#1}}

\newcommand{\N}{\ensuremath{\mathbb{N}}}
\newcommand{\Z}{\ensuremath{\mathbb{Z}}}
\newcommand{\R}{\ensuremath{\mathbb{R}}}

\newcommand{\F}{\ensuremath{\mathcal{F}}}

\newcommand{\Sc}{\ensuremath{\mathbb{S}}}

\newcommand{\ka}{\ensuremath{\rho}}
\newcommand{\kb}{\ensuremath{\sigma}}

\def\map#1#2#3{\mbox{${#1}\!:{#2}\longrightarrow{#3}$}}

\newcommand{\sist}[2]{
  \left\{\!
   \begin{array}{l}
    \dot x=#1 \\[2pt] \dot y=#2
   \end{array}
  \right.
}

\newcommand{\bb}{\ensuremath{b}}

\newcommand{\op}{\ensuremath{\mbox{\rm o}}}

\newcommand{\dsp}{\displaystyle}
\newcommand{\gorro}{\hat}

\title{\bf On the cyclicity of Kolmogorov polycycles
\footnotetext{This work has been partially funded by the Ministry of Science, Innovation and Universities of Spain through the grants PGC2018-095998-B-I00 and MTM2017-86795-C3-2-P and by the Agency for Management of University and Research Grants of Catalonia through the grants 2017SGR1725 and 2017SGR1617.
}}

\author{D. Mar\'{\i}n and J. Villadelprat
\\*[.1truecm]
{\small \textsl{Departament de Matem{\`a}tiques, Edifici Cc,
Universitat Aut{\`o}noma de Barcelona,}}\\*[-.05truecm]
{\small\textsl{08193 Cerdanyola del Vall\`es (Barcelona), Spain}}
\\*[-.05truecm]
{\small \textsl{Centre de Recerca Matem\`atica, Edifici Cc, Campus de Bellaterra,}}\\*[-.05truecm]
{\small \textsl{08193 Cerdanyola del Vall\`es (Barcelona), Spain}}
\\*[.1truecm]
{\small \textsl{Departament d'Enginyeria Inform{\`a}tica i Matem{\`a}tiques, ETSE,}}
\\*[-.05truecm]
{\small \textsl{Universitat Rovira i Virgili, 43007 Tarragona, Spain}}}
                  
\date{\today}

\begin{document}

\maketitle

\begin{abstract}
In this paper we study planar polynomial Kolmogorov's differential systems 
\[
 X_\mu\quad\sist{xf(x,y;\mu),}{yg(x,y;\mu),}
\]
with the parameter $\mu$ varying in an open subset $\Lambda\subset\R^N$. Compactifying $X_\mu$ to the Poincaré disc, the boundary of the first quadrant is an invariant triangle $\Gamma$, that we assume to be a hyperbolic polycycle with exactly three saddle points at its vertices for all $\mu\in\Lambda.$ We are interested in the cyclicity of $\Gamma$ inside the family $\{X_\mu\}_{\mu\in\Lambda},$ i.e., the number of limit cycles that bifurcate from $\Gamma$ as we perturb 
$\mu.$ In our main result we define three functions that play the same role for the cyclicity of the polycycle as the first three Lyapunov quantities for the cyclicity of a focus. As an application we study two cubic Kolmogorov families, with $N=3$ and $N=5$, and in both cases we are able to determine the cyclicity of the polycycle for all $\mu\in\Lambda,$ including those parameters for which the return map along $\Gamma$ is the identity. 
\end{abstract}

\section{Introduction and main results}

The present paper is motivated by the results obtained by Gasull, Mañosa and Mañosas~\cite{GMM02} with regard to the \emph{stability} of an unbounded polycycle $\Gamma$ in the Kolmogorov's polynomial differential systems
\[
 \sist{xf(x,y),}{yg(x,y).}
\]
These systems are widely used in ecology to describe the interaction between two populations, see \cite{Murray} for instance. That being said, the stability of the polycycle is not the main issue to which this paper is addressed. Indeed, assuming that the coefficients of the polynomials $f$ and $g$ depend analytically on a parameter $\mu$, we are interested in the \emph{cyclicity} of the polycycle (see \defic{defcic} below), which roughly speaking is the number of limit cycles that can bifurcate from $\Gamma$ as we perturb $\mu$. In our main result (\teoc{thmA}) we define three functions, $d_1(\mu)$, $d_2(\mu)$ and $d_3(\mu)$, that play the same role for the cyclicity of the polycycle as the first three \emph{Lyapunov quantities} for the cyclicity of a focus. Recall that the displacement map can be analytically extended to a focus and that the Lyapunov quantities are the coefficients of its Taylor's series. On the contrary the displacement map has no smooth extension to a polycycle. At best one can hope that it has some asymptotic expansion. This is indeed the case for the polycycle that we study in the present paper and in order to obtain it we strongly rely in our previous results~\cite{MV19,MV20,MV21} about the asymptotic expansion of the Dulac map of an unfolding of hyperbolic saddles. The principal part of the asymptotic expansion of the displacement map is given in a monomial scale containing a deformation of the logarithm, the so-called Ecalle-Roussarie compensator, and the remainder is uniformly flat with respect to the parameters. The functions $d_i(\mu)$ in \teoc{thmA} are essentially the coefficients of the first three monomials in the principal part, which explains their relation with the cyclicity. For other results regarding the cyclicity of polycycles and more general limit periodic sets the reader is referred to \cite{MDR11,DRR1,Huzak22,RSZ}  and references therein.

Most of the work on planar polynomial differential systems, including this paper, is related to the questions surrounding Hilbert's 16th problem (see for instance \cite{Ily,Roussarie,Ye} and references therein) and its various weakened versions. In this setting it is worth to mention that, using a compactness argument, Roussarie~\cite{Roussarie88} showed that to prove the existential part of Hilbert's 16 problem in the family $\mathcal P_n$ of all polynomial vector fields of degree $\leqslant n$ it is sufficient to show that each \emph{limit periodic set} in $\mathcal P_n$ has finite cyclicity. 

\begin{defi}\label{persistent}
Let $X$ be a vector field on $\R^2$ (or $\Sc^2$). A \emph{graphic} $\Gamma$ for $X$ is a compact, non-empty invariant subset which is a continuous image of $\Sc^1$ and consists of a finite number of isolated singular points $\{p_1,\ldots,p_m,p_{m+1}=p_1\}$ (not necessarily distinct) and compatibly oriented separatrices $\{s_1,\ldots,s_m\}$ connecting them (i.e., such that the $\alpha$-limit set of~$s_j$ is~$p_j$ and the $\omega$-limit set of~$s_j$ is~$p_{j+1}$). A graphic is said to be \emph{hyperbolic} if all its singular points are hyperbolic saddles. A \emph{polycyle} is a graphic with a return map defined on one of its sides. 
\end{defi}

The polycycle that we aim to study is unbounded. In order to investigate the behaviour of the trajectories of a polynomial vector field~$Y$ near infinity we can consider its Poincaré compactification $p(Y)$, see \cite[\S 5]{ADL} for details, which is an analytically equivalent vector field defined on the sphere $\Sc^2$. The points at infinity of~$\R^2$ are in bijective correspondence with the points of the equator of~$\Sc^2$, that we denote by $\ell_{\infty}$. Furthermore the trajectories of $p(Y)$ in $\Sc^2$ are symmetric with respect to the origin and so it suffices to draw its flow in the closed northern hemisphere only, the so called Poincaré disc. 

\begin{defi}\label{defcic}
Let $\{X_\mu\}_{\mu\in\Lambda}$ be a family of vector fields on $\Sc^2$ and suppose that $\Gamma$ is a polycyle for $X_{\mu_0}$. We say that $\Gamma$ has finite cyclicity in the family $\{X_\mu\}_{\mu\in\Lambda}$ if there exist $\kappa\in\N$, $\varepsilon>0$ and $\delta>0$ such that any~$X_\mu$ with $\|\mu-\mu_0\|<\delta$ has at most $\kappa$ limit cycles $\gamma_i$ with $\text{dist}_H(\Gamma,\gamma_i)<\varepsilon.$ The minimum of such $\kappa$ when $\varepsilon$ and $\delta$ go to zero is called the \emph{cyclicity} of $\Gamma$ in $\{X_\mu\}_{\mu\in\Lambda}$ and denoted by $\mathrm{Cycl}\big((\Gamma,X_{\mu_0}),X_\mu\big).$
\end{defi}

\begin{figure}[t]
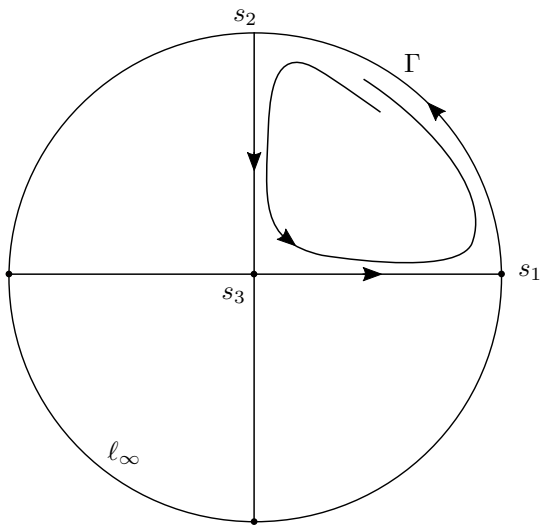

   \centering
  \begin{lpic}[l(0mm),r(0mm),t(0mm),b(5mm)]{dib3(0.75)}
    \lbl[l]{18,13;$\ell_\infty$}   
    \lbl[l]{38,41;$s_3$}
    \lbl[l]{40,90;$s_2$}   
    \lbl[l]{90,45;$s_1$}     
    \lbl[l]{70,82;$\Gamma$}                      
   \end{lpic}
  \caption{Placement of the hyperbolic saddles and the polycycle $\Gamma$ in the Poincaré disc.}
  \label{dib3}
 \end{figure}

In this paper we consider the family of vector fields $\{X_\mu\}_{\mu\in\Lambda}$ given by
\begin{equation}\label{kol}
 X_\mu\!:=f(x,y;\mu)x\partial_x+g(x,y;\mu)y\partial_y
\end{equation}
where $\Lambda$ is an open subset of $\R^N$ and $f$ and $g$ are polynomials in $x$ and $y$ of degree $n\in\N$ with the coefficients depending analytically on $\mu$. The standing hypothesis on the family $\{X_\mu\}_{\mu\in\Lambda}$ are the following:
\begin{enumerate}

\item[\bf{H1}] $f(z,0;\mu)>0$, $g(0,z;\mu)<0$ and $\big(f_n-g_n\big)(1,z;\mu)<0$ for all $z>0$ and $\mu\in\Lambda.$

\item[\bf{H2}] $\lambda_1(\mu)\!:=\left(\frac{f_n}{g_n-f_n}\right)(1,0;\mu)$, $\lambda_2(\mu)\!:=\left(\frac{f_n-g_n}{g_n}\right)(0,1;\mu)$ and $\lambda_3(\mu)\!:=-\left(\frac{g}{f}\right)(0,0;\mu)$,  are well defined and strictly positive for all $\mu\in \Lambda.$

\end{enumerate}
Here, and in what follows, $f_n$ and $g_n$ denote, respectively, the homogeneous part of degree $n$ of $f$ and $g.$
Conditions {\bf H1} and {\bf H2} guarantee that, after compactifying the polynomial vector field $X_\mu$ to the Poincaré disc, the boundary of the first quadrant is a polycyle with three hyperbolic saddles, see \figc{dib3},
\[
 \text{$s_1\!:=\{y=0\}\cap \ell_\infty$, $s_2\!:=\{x=0\}\cap \ell_\infty$ and $s_3\!:=(0,0)$.}
\] 
From now on we shall denote this polycyle by $\Gamma$, that we remark is a compact subset of the Poincaré disc. The hyperbolicity ratios of the saddles at its vertices are precisely the ones given in ${\bf H2}$.
We also define:

\begin{equation*}
\begin{array}{ll}
 \dsp L_{11}(u)=\exp\left(\int_0^u\left(\left(\frac{f-g}{f}\right)(1/z,0)+\frac{1}{\lambda_1}\right)\frac{dz}{z}\right)
 &
 \dsp L_{12}(u)=\exp\left(\int_0^u\left(\left(\frac{f_n}{f_n-g_n}\right)(1,z)+\lambda_1\right)\frac{dz}{z}\right)
 \\[15pt]
 \dsp L_{21}(u)=\exp\left(\int_0^u\left(\left(\frac{g_n}{g_n-f_n}\right)(z,1)+\frac{1}{\lambda_2}\right)\frac{dz}{z}\right)
 &
 \dsp L_{22}(u)=\exp\left(\int_0^u\left(\left(\frac{g-f}{g}\right)(0,1/z)+\lambda_2\right)\frac{dz}{z}\right)
 \\[15pt]
 \dsp L_{31}(u)=\exp\left(\int_0^u\left(\left(\frac{g}{f}\right)(z,0)+\lambda_3\right)\frac{dz}{z}\right)
 &
 \dsp L_{32}(u)=\exp\left(\int_0^u\left(\left(\frac{f}{g}\right)(0,z)+\frac{1}{\lambda_3}\right)\frac{dz}{z}\right)
 \\
 \end{array}
\end{equation*}
together with
\begin{equation}\label{def_M}
\text{$M_1(u)=-\frac{L_{11}(u)}{u}\,\partial_2\!\left(\frac{g}{f}\right)(1/u,0)$ 
and $M_3(u)=L_{31}(u)\,\partial_2\!\left(\frac{g}{f}\right)(u,0)$.}
\end{equation}
We point out that all these functions depend on the parameter $\mu.$ This dependence is omitted for the sake of shortness when there is no risk of confusion. 

We can now state our main result, which is addressed to the cyclicity of the polycycle $\Gamma$ inside the polynomial family $\{X_\mu\}_{\mu\in\Lambda}.$  More formally, we should refer to the compactified family $\{p(X_\mu)\}_{\mu\in\Lambda}$ of vector fields on $\Sc^2$ but for the simplicity in the exposition we commit an abuse of language by identifying both families. It is clear that the number of limit cycles of $p(X_\mu)$ and $X_\mu$ is the same. In the statement $\mathscr R(\,\cdot\,;\mu)$ stands for the return map of the vector field~$X_\mu$ around the polycycle $\Gamma$ (see \figc{dib3}) and we use the notion of functional independence that is given in \defic{indepe}.

\begin{bigtheo}\label{thmA}
Let us consider the family of Kolmogorov polynomial vector fields $\{X_\mu\}_{\mu\in\Lambda}$ given in \refc{kol} and verifying the assumptions {\bf H1} and {\bf H2}. Then, for any $\mu_0\in\Lambda,$ the following assertions hold with regard to the cyclicity of the polycycle $\Gamma$ inside the family:

\begin{enumerate}[$(a)$]
\item If $d_1(\mu)\!:=1-\lambda_1\lambda_2\lambda_3$ does not vanish at $\mu_0$ then 
         $\mathrm{Cycl}\big((\Gamma,X_{\mu_0}),X_\mu\big)=0.$
         
\item If $d_1$ vanishes and is independent at $\mu_0$ and $\mathscr R(\,\cdot\,;\mu_0)\not\equiv \text{Id}$ then        
         $\mathrm{Cycl}\big((\Gamma,X_{\mu_0}),X_\mu\big)\geqslant 1.$

\item If $d_2(\mu)\!:=\log\left(\left(\frac{L_{12}}{L_{21}}\right)^{\lambda_2}
         \left(\frac{L_{31}}{L_{11}}\right)^{\lambda_1\lambda_2}\frac{L_{22}}{L_{32}}\right)(1)$
         does not vanish at $\mu_0$ 
         then $\mathrm{Cycl}\big((\Gamma,X_{\mu_0}),X_\mu\big)\leqslant 1.$

\item If $d_1$ and $d_2$ vanish and are independent at $\mu_0$ and $\mathscr R(\,\cdot\,;\mu_0)\not\equiv \text{Id}$ then 
         $\mathrm{Cycl}\big((\Gamma,X_{\mu_0}),X_\mu\big)\geqslant 2.$
\end{enumerate}
In case that $\lambda_1(\mu_0)<1$, $\lambda_2(\mu_0)>1$ and $\lambda_3(\mu_0)>1$ 
then the following assertions hold as well:
\begin{enumerate}[$(a)$]\setcounter{enumi}{4}
\item If $d_3(\mu)\!:=\hat M_3(\lambda_3,1)L_{11}(1)-\hat M_1\big(\frac{1}{\lambda_1},1\big)L_{31}(1)$ 
         does not vanish at $\mu_0$ then 
         $\mathrm{Cycl}\big((\Gamma,X_{\mu_0}),X_\mu\big)\leqslant 2.$
         
\item If $d_1,$ $d_2$ and $d_3$ vanish and are independent at $\mu_0$ and 
         $\mathscr R(\,\cdot\,;\mu_0)\not\equiv \text{Id}$ then
         $\mathrm{Cycl}\big((\Gamma,X_{\mu_0}),X_\mu\big)\geqslant 3.$
\end{enumerate}
\end{bigtheo}


Let us make some remarks with regard to the regularity of the functions $d_1$, $d_2$ and $d_3$ defined in the statement. On account of the hypothesis {\bf H1} and {\bf H2} it is evident that $d_1$ is analytic on the whole parameter space~$\Lambda$. On the other hand, $d_2$ is defined in terms of the functions~$\mu\mapsto L_{ij}(1)$, which in turn are given by some (apparently) improper integrals. By applying the Weierstrass Division Theorem one can easily show that each $L_{ij}(1)$ is an analytic strictly positive function, so that~$d_2$ is also analytic on~$\Lambda.$ Finally, $d_3$ is given by means of a sort of incomplete Mellin transform (which is defined in \propc{L8}) of the functions~$M_1$ and~$M_3$ in~\refc{def_M}. One can  show that the hypothesis $\mathbf{H1}$ and $\mathbf{H2}$ imply that each $M_i(u;\mu)$ is analytic on $(-\varepsilon,+\infty)\times\Lambda$ for some $\varepsilon>0.$ Taking this into account, by applying~$(d)$ in \propc{L8} it follows that $d_3$ is a meromorphic function on $\Lambda$ having poles only at those $\mu_0$ such that $1/\lambda_1(\mu_0)\in\N$ or $\lambda_3(\mu_0)\in\N$.

Also with regard to the statement of \teoc{thmA}, the assertions $(e)$ and $(f)$ hold under the assumptions $\lambda_1(\mu_0)<1$, $\lambda_2(\mu_0)>1$ and $\lambda_3(\mu_0)>1$. However one can always reduce to this case provided that $\lambda_i(\mu_0)\neq 1$ for $i=1,2,3$ by means of a rescaling of time and a projective change of coordinates that permute conveniently the three singular points of the polycycle. 

The paper \cite{GMM02} constitutes an important previous contribution to the study of Kolmogorov polycycles that should be referred. Indeed, following our notations and definitions, the authors prove (see  \cite[Theorem~1]{GMM02}) that if $d_1(\mu_0)=0$ then the return map of $X_{\mu_0}$ around the polycycle $\Gamma$ is of the form 
\begin{equation}\label{v1}
 \mathscr R(s;\mu_0)=\Delta s+\op(s),
\end{equation}
cf. $(b)$ in \teoc{prop1}, and they also provide the explicit expression of the coefficient $\Delta.$ This coefficient is given as the limit  of a sum of three improper integrals, which computed separately diverge. An easy manipulation of the integrals shows that these divergences cancel each other, yielding to the expression of~$d_2$ given in \teoc{thmA}. It is important to remark that the expansion in \refc{v1} can not be used to obtain an upper bound for $\mathrm{Cycl}\big((\Gamma,X_{\mu_0}),X_\mu\big)$ because the remainder is not uniform with respect to the parameters. 
It is possible however to use it to obtain lower bounds. In this direction the authors prove in \cite[Corollary~5]{GMM02}
that if $d_1$ vanishes and is independent at $\mu_0$ and $d_2(\mu_0)\neq 0$ then $\mathrm{Cycl}\big((\Gamma,X_{\mu_0}),X_\mu\big)\geqslant1$. Since $d_2(\mu_0)\neq 0$ implies $\mathscr R(\,\cdot\,;\mu_0)\not\equiv Id$ by \teoc{prop1}, this lower bound follows by applying $(b)$ in \teoc{thmA}.

The paper is organised in the following way. \secc{provateoA} is entirely devoted to prove \teoc{thmA} and for that purpose we rely in our previous results about the asymptotic expansion of the Dulac map of a hyperbolic saddle that we obtain in \cite{MV19,MV20,MV21}. For this reason, before starting the proof of \teoc{thmA} we first state these results and introduce the necessary definitions. 
The asymptotic expansion of the displacement map near the polycycle is given in \teoc{prop1} and constitutes the fundamental tool in order to prove \teoc{thmA}. As a by-product of this expansion we obtain a method to study the stability of the polycycle, see \obsc{rem_stab}. \secc{appli} is addressed to the applications. The first one is \teoc{ex1}, where we consider a Kolmogorov's cubic system depending on three parameters that was previously studied in~\cite{GMM02}. The authors in that paper show that there exist parameters for which the cyclicity of the polycycle is at least 1. In the present paper we obtain the exact cyclicity of all the parameters in the family (that can be 0 or 1), including the case in which the return map along the polycycle is the identity. We also show that there exists exactly one singularity in the first quadrant, which can be a focus or a center, and we compute its cyclicity. Finally we prove that it is not possible a simultaneous bifurcation of limit cycles from the polycycle and that singularity. We give our second application in \teoc{ex2}, where we consider a Kolmogorov's cubic system depending on five parameters. In this case we also provide the exact cyclicity of $\Gamma$ for all the parameters in the family, which again can be 0 or 1. 

\section{Proof of \teoc{thmA}}\label{provateoA}

In order to tackle the proof of \teoc{thmA} we will appeal to some previous results from \cite{MV19,MV20,MV21} about the asymptotic expansion of the Dulac map. For reader's convenience we gather these results in \propc{3punts}. To this end it is first necessary to introduce some new notation and definitions. 

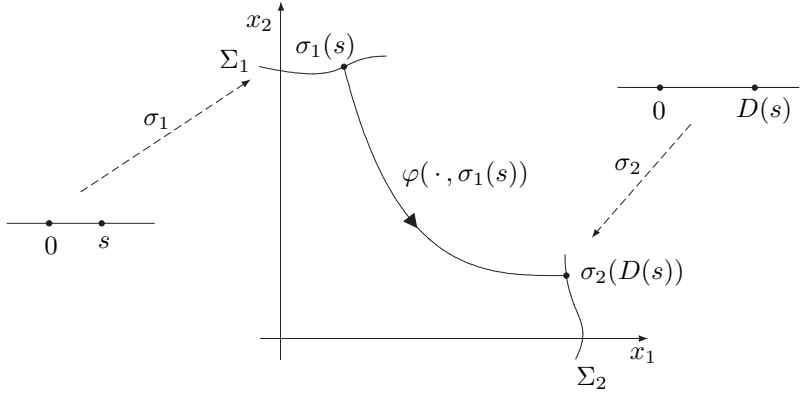
\begin{figure}[t]
   \centering
  \begin{lpic}[l(0mm),r(0mm),t(0mm),b(5mm)]{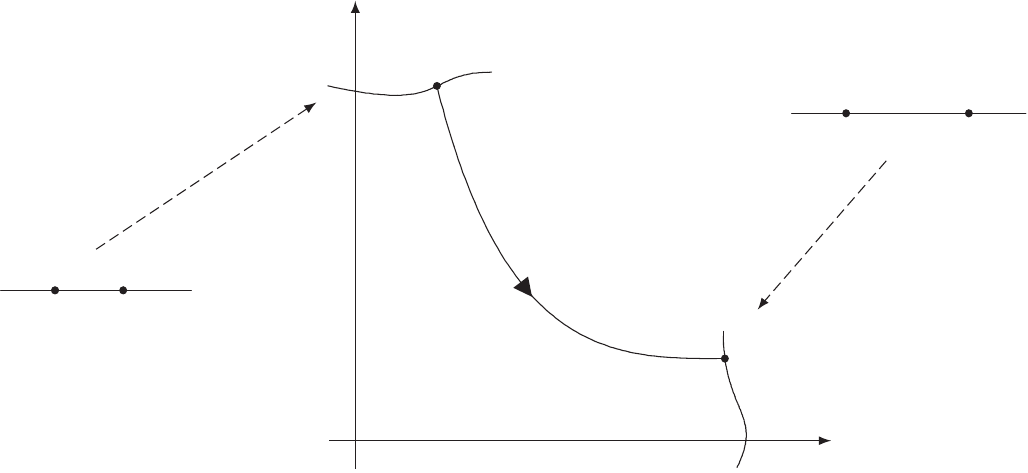}
   \lbl[l]{5,15.5;$0$}   
   \lbl[l]{12,16;$s$}   
   \lbl[l]{18,32;$\sigma_1$}
   \lbl[l]{28,40;$\Sigma_1$}      
   \lbl[l]{31.5,45;$x_2$}   
   \lbl[l]{38,42;$\sigma_1(s)$}   
   \lbl[l]{52,25;$\varphi(\,\cdot\,,\sigma_1(s))$}
   \lbl[l]{75.5,12;$\sigma_2(D(s))$}  
   \lbl[l]{75,-2;$\Sigma_2$}   
   \lbl[l]{82,1;$x_1$}
   \lbl[l]{80,26;$\sigma_2$}
   \lbl[l]{85,33.5;$0$}
   \lbl[l]{96,33.5;$D(s)$}                   
   \end{lpic}
  \caption{Definition of the Dulac map $D(\,\cdot\,;{\hat\nu})$, where $\varphi(t,p;{\hat\nu})$ is the solution of~$X_{\hat\nu}$ passing through the point $p\in\mathscr U$ at time $t=0.$}
  \label{DefTyR}
 \end{figure}

Setting $\hat\nu\!:=(\lambda,\nu)\in\hat W\!:=(0,+\infty)\times W$ with $W$ an open set of $\R^N,$ we consider the family of vector fields $\{X_{\hat\nu}\}_{{\hat\nu}\in\hat W}$ with
\begin{equation}\label{X}
 X_{\hat\nu}({x_1},{x_2})={x_1}P_1({x_1},{x_2};{\hat\nu})\partial_{x_1}+{x_2}P_2({x_1},{x_2};{\hat\nu})\partial_{x_2}
\end{equation}
where
\begin{itemize}
\item $P_1$ and $P_2$ belong to $\mathscr C^{\omega}(\mathscr U\!\times\!\hat W)$ for some open set $\mathscr U$ of $\R^2$ containing the origin, 
\item $P_1({x_1},0;{\hat\nu})>0$ and $P_2(0,{x_2};{\hat\nu})<0$ for all $({x_1},0),(0,{x_2})\in\mathscr U$ and ${\hat\nu}\in \hat W,$
\item $\lambda=-\frac{P_2(0,0;\nu)}{P_1(0,0;\nu)}$.
\end{itemize}
Thus, for all $\hat\nu\in\hat W$, the origin is a hyperbolic saddle of $X_{\hat\nu}$ with the separatrices lying in the axis. We point out that here the hyperbolicity ratio of the saddle is an independent parameter, although in the proof of \teoc{thmA} we will have $\lambda=\lambda(\nu)$. The reason for this is that the hyperbolicity ratio turns out to be the ruling parameter in our results and, besides, having it uncoupled from the rest of parameters simplifies the notation in the statements.
Moreover, for $i=1,2,$ we consider a $\mathscr C^{\omega}$~transverse section \map{\sigma_i}{(-\varepsilon,\varepsilon)\times \hat W}{\Sigma_i} to~$X_{{\hat\nu}}$ at $x_i=0$ defined by
 \[
  \sigma_i(s;{\hat\nu})=\bigl(\sigma_{i1}(s;{\hat\nu}),\sigma_{i2}(s;{\hat\nu})\bigr)
 \]
such that $\sigma_1(0,{\hat\nu})\in\{(0,x_2);x_2>0\}$ and $\sigma_2(0,{\hat\nu})\in\{(x_1,0);x_1>0\}$
for all ${\hat\nu}\in \hat W.$ We denote the Dulac map of~$X_{\hat\nu}$ from $\Sigma_1$ to $\Sigma_2$ by $D(\,\cdot\,;{\hat\nu})$, see \figc{DefTyR}. 
The asymptotic expansion of $D(s;\hat\nu)$ at $s=0$ consists of a remainder and a principal part. The principal part is given in a monomial scale that contains a deformation of the logarithm, the so-called Ecalle-Roussarie compensator, whereas the remainder has good flatness properties with respect to the parameters. We next give precise definitions of these key notions. 

\begin{defi}\label{defi_comp}
The function defined for $s>0$ and $\alpha\in\R$ by means of
 \[
  \omega(s;\alpha)\:=
  \left\{
   \begin{array}{ll}
    \frac{s^{-\alpha}-1}{\alpha} & \text{if $\alpha\neq 0,$}\\[2pt]
    -\log s & \text{if $\alpha=0,$}
   \end{array}
  \right.
 \]
is called the \emph{Ecalle-Roussarie compensator}. 
\end{defi}

\begin{defi}
Consider an open subset $U\subset\hat W\subset \R^{N+1}.$ We say that a function $\psi(s;{\hat\nu})$ belongs to the class $\mathscr C^\infty_{s>0}(U)$ if there 
exist an open neighbourhood $\Omega$ of 
\[ 
 \{(s,{\hat\nu})\in\R^{N+2};s=0,{\hat\nu}\in U\}=\{0\}\times U
\]  
in $\R^{N+2}$ such that $(s,{\hat\nu})\mapsto \psi(s;{\hat\nu})$ is $\mathscr C^\infty$ on $\Omega\cap\big((0,+\infty)\times U\big)$. 
\end{defi}

More formally, the definition of $\mathscr C^\infty_{s>0}(U)$ must be thought in terms of germs with respect to relative neighbourhoods of $\{0\}\times U$ in $(0,+\infty)\times U$. In doing so $\mathscr C^\infty_{s>0}(U)$ becomes a ring. 

We can now introduce the notion of flatness that we shall use in the sequel. 

\begin{defi}\label{defi2} 
Consider an open subset $U\subset\hat W\subset\R^{N+1}.$ Given $L\in\R$ and ${\hat\nu}_0\in U$, we say that a function $\psi(s;{\hat\nu})\in\mathscr C^\infty_{s>0}(U)$ is \emph{$L$-flat with respect to $s$ at ${\hat\nu}_0$}, and we write $\psi\in\F_L^\infty({\hat\nu}_0)$, if for each ${\hat\nu}=({\hat\nu}_0,\ldots,{\hat\nu}_{N+1})\in\Z_{\geq 0}^{N+2}$  with $|{\hat\nu}|={\hat\nu}_0+\cdots+{\hat\nu}_{N+1}\leqslant K$ there exist a neighbourhood~$V$ of ${\hat\nu}_0$ and $C,s_0>0$ such that
\begin{equation*}
 \left|\frac{\partial^{|{\hat\nu}|}\psi(s;{\hat\nu})}{\partial s^{{\hat\nu}_0}
 \partial{\hat\nu}_1^{{\hat\nu}_1}\cdots\partial{\hat\nu}_{N+1}^{{\hat\nu}_{N+1}}}\right|\leqslant C s^{L-{\hat\nu}_0}
 \text{ for all $s\in(0,s_0)$ and ${\hat\nu}\in V$.} 
\end{equation*}
If $W$ is a (not necessarily open) subset of $U$ then define $\F_L^\infty(W)\!:=\bigcap_{{\hat\nu}_0\in W}\F_L^\infty({\hat\nu}_0).$
\end{defi}

Apart from the remainder and the monomial order, the most important ingredient for our purposes is the explicit expression of the coefficients in the asymptotic expansion. In order to give them we introduce next some additional notation, where for the sake of shortness the dependence on ${\hat\nu}=(\lambda,\nu)$ is omitted. We define the functions:
\begin{equation}\label{def_fun}
\begin{array}{ll}
\dsp L_1(u)\!:=\exp\int_0^u\left(\frac{P_1(0,z)}{P_2(0,z)}+\frac{1}{\lambda}\right)\frac{dz}{z} & 
\dsp L_2(u)\!:=\exp\int_0^u\left(\frac{P_2(z,0)}{P_1(z,0)}+{\lambda}\right)\frac{dz}{z} \\[15pt]
\dsp M_1(u)\!:=L_1(u)\partial_1\!\left(\frac{P_1}{P_2}\right)(0,u)&  
\dsp M_2(u)\!:=L_2(u)\partial_2\!\left(\frac{P_2}{P_1}\right)(u,0)\\[15pt]
\end{array}
\end{equation}
On the other hand, for shortness as well, we use the compact notation $\sigma_{ijk}$ for the $k$th derivative at $s=0$ of the $j$th component of $\sigma_i(s;{\hat\nu})$, i.e., 
\[
 \sigma_{ijk}({\hat\nu})\!:=\partial^k_s\sigma_{ij}(0;{\hat\nu}).
\]
Taking this notation into account we also introduce the following real values, where once again we omit the dependence on ${\hat\nu}$:
\begin{equation}\label{def_S}
\begin{array}{l}
\dsp S_1\!:
=\frac{\sigma_{112}}{2\sigma_{111}}-\frac{\sigma_{121}}{\sigma_{120}}\left(\frac{P_1}{P_2}\right)\!(0,\sigma_{120})-\frac{\sigma_{111}}{L_1(\sigma_{120})}\gorro{M}_1(1/\lambda,\sigma_{120})\\[20pt]
\dsp S_2\!:
=\frac{\sigma_{222}}{2\sigma_{221}}-\frac{\sigma_{211}}{\sigma_{210}}\left(\frac{P_2}{P_1}\right)\!(\sigma_{210},0)-{\frac{\sigma_{221}}{L_2(\sigma_{210})}}\gorro{M}_2(\lambda,\sigma_{210}).
\end{array}
\end{equation}
Here $\hat M_i$ stands for a sort of incomplete Mellin transform of $M_i$ that will be defined by \propc{L8} below. 
We can now state the following result, which gathers Theorem~A and Theorem~4.1 in~\cite{MV21} and that it will constitute the key tool in order to prove the main result in the present paper.

\begin{prop}\label{3punts}
Let $D(s;{\hat\nu})$ be the Dulac map of the hyperbolic saddle \refc{X} from $\Sigma_1$ and $\Sigma_2$ and define
\[
\Delta_{0}(\hat\nu)=\frac{\sigma_{111}^\lambda\sigma_{120}}{L_1^\lambda(\sigma_{120})}\frac{L_2(\sigma_{210})}{\sigma_{221}\sigma_{210}^\lambda},\text{ }
\Delta_{1}(\hat\nu)=\Delta_{0}\lambda S_1\text{ and } 
\Delta_{2}(\hat\nu) =-\Delta_{0}^2S_2,
\]
where $\hat\nu=(\lambda,\nu)\in\hat W=(0,+\infty)\times W.$
Then $\Delta_{0}$ is analytic and strictly positive on $\hat W$, 
$\Delta_{1}$ is meromorphic on $\hat W$ with poles only at $\lambda\in\frac{1}{\N}$ and 
$\Delta_{2}$ is meromorphic on $\hat W$ with poles only at $\lambda\in\N.$
Moreover
the following assertions hold:
\begin{enumerate}[$(1)$]

\item If $\lambda_0<1$ then $D(s;{\hat\nu})=s^\lambda\big(\Delta_{0}({\hat\nu}) +\Delta_{2}({\hat\nu})s^{\lambda}+\F_{\ell}^\infty(\{\lambda_0\}\times W)\big)$ for any ${\ell}\in  \big[\lambda_0,\min(2\lambda_0,1)\big)$.

\item If $\lambda_0=1$ then $D(s;{\hat\nu})=s^\lambda\big(\Delta_{0}({\hat\nu}) +\boldsymbol{\Delta}^{\lambda_0}(\omega;{\hat\nu})s +\F_{\ell}^\infty(\{\lambda_0\}\times W)\big)$ for any ${\ell}\in [1,2),$ where 
\[
 \boldsymbol{\Delta}^{\lambda_0}(\omega;{\hat\nu})=\Delta_{1}({\hat\nu})+\Delta_{2}({\hat\nu})(1+\alpha\omega),
\] 
$\alpha=1-\lambda$ and $\omega=\omega(s;\alpha)$.

\item If $\lambda_0>1$ then $D(s;{\hat\nu})=s^\lambda\big(\Delta_{0}({\hat\nu}) +\Delta_{1}({\hat\nu})s +\F_{\ell}^\infty(\{\lambda_0\}\times W)\big)$ for any ${\ell}\in \big[1,\min(\lambda_0,2)\big)$.

\end{enumerate}
In particular $D(s;{\hat\nu})=s^\lambda\big(\Delta_{0}({\hat\nu}) +\F_{\ell}^\infty(\{\lambda_0\}\times W)\big)$ for any ${\ell}\in \big(0,\min(\lambda_0,1)\big)$.
\end{prop}

The flatness $\ell$ of the remainder can range in a certain interval depending on $\lambda_0.$ The left endpoint of this interval is only given for completeness to guarantee that all the monomials in the principal part are relevant (i.e., they cannot be included in the remainder). The important information about the flatness is given by the right endpoint.
A key tool in order to give a closed expression of the coefficients $\Delta_i$ is the use of a sort of incomplete Mellin transform, which is accurately defined in the next result. For a proof of this result the reader is referred to \cite[Appendix B]{MV21}.

\begin{prop}\label{L8}
Let us consider an open interval $I$ of $\R$ containing $x=0$ and an open subset $U$ of $\R^M$.
\begin{enumerate}[$(a)$]

\item Given $f(x;{\upsilon})\in\mathscr C^{\infty}(I\times U)$, there exits a unique $\hat f(\alpha,x;{\upsilon})\in\mathscr C^{\infty}((\R\setminus\Z_{\ge 0})\times I\times U)$ such that 
\begin{equation*}
 x\partial_x\hat f({\alpha},x;{\upsilon})-\alpha\hat f({\alpha},x;{\upsilon})=f(x;{\upsilon}).
\end{equation*}

\item If $x\in I\setminus\{0\}$ then $\partial_x(\hat f({\alpha},x;{\upsilon})|x|^{-\alpha})=f(x;{\upsilon})\frac{|x|^{-\alpha}}{x}$ and, taking any $k\in\Z_{\ge0}$ with $k>\alpha$,
\begin{equation*}
\hat f(\alpha,x;{\upsilon})=
\sum_{i=0}^{k-1}\frac{\partial_x^if(0;{\upsilon})}{i!(i-\alpha)}x^i+|x|^{\alpha}\int_0^x\!\left(f(s;{\upsilon})-T_0^{k-1}f(s;{\upsilon})\right)|s|^{-\alpha}\frac{ds}{s},
\end{equation*}
where $T_0^kf(x;{\upsilon})=\sum_{i=0}^{k}\frac{1}{i!}\partial_x^if(0;{\upsilon})x^i$ is the $k$-th degree Taylor polynomial of $f(x;{\upsilon})$ at $x=0$.

\item 
For each $(i_0,x_0,{\upsilon}_0)\in\Z_{\ge 0}\times I\times W$ the function $(\alpha,x,{\upsilon})\mapsto(i_0-\alpha)\hat f(\alpha,x;{\upsilon})$ extends $\mathscr C^\infty$ at $(i_0,x_0,{\upsilon}_0)$ and, moreover, it tends to $\frac{1}{i_0!}\partial_x^{i_0}f(0;{\upsilon}_0)x_0^{i_0}$ as $(\alpha,x,{\upsilon})\to (i_0,x_0,{\upsilon}_0).$

\item If $f(x;{\upsilon})$ is analytic on $I\times U$ then $\hat f(\alpha,x;{\upsilon})$ is analytic on $(\R\setminus\Z_{\ge 0})\times I\times U$. Finally,
for each $(\alpha_0,x_0,{\upsilon}_0)\in\Z_{\ge 0}\times I\times U$ the function
$(\alpha,x,{\upsilon})\mapsto(\alpha_0-\alpha)\hat f(\alpha,x;{\upsilon})$ extends analytically to $(\alpha_0,x_0,{\upsilon}_0)$.
\end{enumerate}
\end{prop}

On account of this result for each $M_i(u;\hat\nu)$ in \refc{def_fun} we have that $(\alpha,u;\hat\nu)\mapsto\hat M_i(\alpha,u;\hat\nu)$ is a well defined meromorphic function with poles only at $\alpha\in\Z_{\geq 0}$. Accordingly, see \refc{def_S}, $\hat M_1(1/\lambda,\sigma_{120})$ and $\hat M_2(\lambda,\sigma_{210})$ are the values (depending on $\hat\nu$) that we obtain by taking $\hat M_1(\alpha,u;\hat\nu)$ with $\alpha=1/\lambda$ and $u=\sigma_{120}(\hat\nu)$ and by taking $\hat M_2(\alpha,u;\hat\nu)$ with $\alpha=\lambda$ and $u=\sigma_{210}(\hat\nu),$ respectively.

\begin{figure}[t]
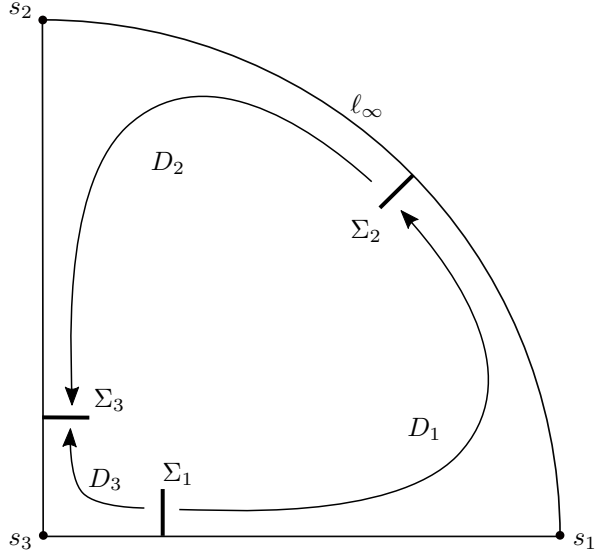

   \centering
  \begin{lpic}[l(0mm),r(0mm),t(0mm),b(5mm)]{dib2(0.75)}
    \lbl[l]{22,12;$\Sigma_1$}
    \lbl[l]{10,25;$\Sigma_3$}   
    \lbl[l]{55,55;$\Sigma_2$}       
    \lbl[l]{55,77;$\ell_\infty$}   
    \lbl[l]{65,20;$D_1$}
    \lbl[l]{9,11;$D_3$}   
    \lbl[l]{20,67;$D_2$}    
    \lbl[l]{94,0;$s_1$}
    \lbl[l]{-5,94;$s_2$}   
    \lbl[l]{-5,0;$s_3$}                      
   \end{lpic}
  \caption{Auxiliary Dulac maps for the definition of $\mathscr D=D_2\circ D_1-D_3$ in \teoc{prop1}. The return map 
               in \teoc{thmA}, with respect to the transverse section $\Sigma_1$, would be $\mathscr R=D_3^{-1}\circ D_2\circ 	       D_1=D_3^{-1}\circ\mathscr D+\text{Id}$.}
  \label{dib2}
 \end{figure}

At this point we get back to the setting treated in the present paper and from now on we recover the original notation for the parameters in the family under consideration, see \refc{kol}.
In order to study the Dulac maps of the hyperbolic saddles at the vertices of the polycycle $\Gamma$ we take three local transverse sections $\Sigma_1$, $\Sigma_2$, and $\Sigma_3$ parametrised, respectively, by $s\mapsto (1,s)$, $s\mapsto (1/s,1/s)$ and $s\mapsto (s,1)$ with $s>0.$ We define $D_1(s;\mu)$ to be the Dulac map of $X_\mu$ from $\Sigma_1$ to $\Sigma_2,$ $D_2(s;\mu)$ to be the Dulac map of $X_\mu$ from $\Sigma_2$ to $\Sigma_3$ and, finally, $D_3(s;\mu)$ to be the Dulac map of $-X_\mu$ from $\Sigma_1$ to $\Sigma_3,$ see \figc{dib2} 
 x
It is then clear that the limit cycles of~$X_\mu$ near 
$\Gamma$ are in one to one correspondence with the isolated positive zeroes of 
 \[
  \mathscr D(s;\mu)\!:=\big(D_2\circ D_1-D_3\big)(s;\mu)
 \]
near $s=0.$ The proof of \teoc{thmA} strongly relies in our next result, where we get the asymptotic expansion of $\mathscr D(s;\mu)$ at $s=0$ and we compute its coefficients. In its statement $d_i(\mu)$, for $i=1,2,3,$ are the functions defined in \teoc{thmA} and $\R\{\mu\}_{\mu_0}$ stands for the local ring of convergent power series at $\mu_0.$

\begin{theo}\label{prop1} Let us fix any $\mu_0\in\Lambda$ and set $\lambda_i^0\!:=\lambda_i(\mu_0)$ for $i=1,2,3.$

\begin{enumerate}[$(a)$]

\item If $\lambda_1^0\lambda_2^0\lambda_3^0\neq 1$ then, for any $\ell_1\in\big(\min(\lambda_1^0\lambda_2^0,1/\lambda_3^0),\min(\lambda_1^0+\lambda_1^0\lambda_2^0,
1+\lambda_1^0\lambda_2^0,2\lambda_1^0\lambda_2^0,1+1/\lambda_3^0,2/\lambda_3^0)\big)$,
\begin{equation*}
 \mathscr D(s;\mu)=a_1(\mu)s^{\lambda_1\lambda_2}-a_2(\mu)s^{1/\lambda_3}+\F_{\ell_1}^\infty(\mu_0),
\end{equation*}
where $a_1$ and $a_2$ are analytic and strictly positive functions on $\Lambda.$

\item If $\lambda_1^0\lambda_2^0\lambda_3^0=1$ then, for any 
$\ell_2\in(0,\min(1,\lambda_1^0,\lambda_1^0\lambda_2^0))$,
\begin{equation*}
 \mathscr D(s;\mu)=\big(b_1(\mu)\,\omega\big(s;\alpha(\mu)\big)+b_2(\mu)+\F_{\ell_2}^\infty(\mu_0)\big)s^{1/\lambda_3},
\end{equation*}
where $\alpha=1/\lambda_3-\lambda_1\lambda_2$, $b_1=\alpha a_1$ and $b_2=a_1-a_2.$ Moreover the equalities 
\[
 (b_1)=(d_1)\text{ and }(b_1,b_2)=(d_1,d_2)
\]
between ideals over the ring $\R\{\mu\}_{\mu_0}$ are verified. 
Assuming $\lambda_1^0>1$, $\lambda_2^0>1$ and $\lambda_3^0<1$ additionally then, for any $\ell_3\in \big(1,\min(2,\lambda_1^0,1/\lambda_3^0)\big)$,
\begin{equation*}
 \mathscr D(s;\mu)=\big(b_1(\mu)\,\omega\big(s;\alpha(\mu)\big)+b_2(\mu)+b_3(\mu)s+\F_{\ell_3}^\infty(\mu_0)\big)s^{1/\lambda_3}U(s;\mu),
\end{equation*}
where $b_3$ is an analytic function at $\mu_0$ verifying that 
\[
 (b_1,b_2,b_3)=(d_1,d_2,d_3)
\] 
over the ring $\R\{\mu\}_{\mu_0}$ and $U$ is an analytic function such that $U(0;\mu_0)=1$. 
\end{enumerate}
\end{theo}

\begin{prova}
In order to study the Dulac map $D_1$ from $\Sigma_1$ to $\Sigma_2$ we compactify $X_\mu$ by means of the coordinate change $\{x_1=\frac{y}{x},x_2=\frac{1}{x}\}.$ One can easily verify that the new vector field is orbitally conjugated to \refc{X} particularised with $P_1(x_1,x_2)=x_2^n\big(f-g)(\frac{1}{x_2},\frac{x_1}{x_2})$ and $P_2(x_1,x_2)=x_2^nf(\frac{1}{x_2},\frac{x_1}{x_2}),$ whereas in these coordinates the transverse sections $\Sigma_1$ and $\Sigma_2$ are parametrised by $\sigma_1(s)=(s,1)$ and $\sigma_2(s)=(1,s)$, respectively. The hyperbolicity ratio of the saddle at the origin is $\lambda_1=-\big(\frac{P_2}{P_1}\big)(0,0)=\big(\frac{f_n}{g_n-f_n}\big)(1,0).$ Therefore, by applying \propc{3punts} we can assert that
\begin{equation}\label{prop1eq1}
 D_1(s)=\Delta_{10}s^{\lambda_1}\big(1+\F^\infty_{{\ell}_1}(\mu_0)\big),\text{ with $\Delta_{10}\!:=\big(L_{12}L_{11}^{-\lambda_1}\big)(1)$}
\end{equation}
and any $0<{\ell}_1<\min(1,\lambda_1^0)$. Recall here that $\lambda_i^0=\lambda_i(\mu_0)$ for $i=1,2,3$ by definition.

Next, to analyse the Dulac map $D_2$ from $\Sigma_2$ to $\Sigma_3$ we compactify $X_\mu$ performing the change of coordinates given by $\{x_1=\frac{1}{y},x_2=\frac{x}{y}\}.$ One can check that the new vector field is orbitally conjugated to \refc{X}  with $P_1(x_1,x_2)=x_1^ng(\frac{x_2}{x_1},\frac{1}{x_1})$ and
$P_2(x_1,x_2)=x_1^n\big(g-f)(\frac{x_2}{x_1},\frac{1}{x_1})$ and that in these coordinates the transverse sections $\Sigma_2$ and $\Sigma_3$ are parametrised by $\sigma_1(s)=(s,1)$ and $\sigma_2(s)=(1,s)$, respectively. The hyperbolicity ratio of the saddle at the origin is $\lambda_2=-\big(\frac{P_2}{P_1}\big)(0,0)=\big(\frac{f_n-g_n}{g_n}\big)(0,1).$ Thus, by \propc{3punts} again,
\begin{equation}\label{prop1eq2}
 D_2(s)=\Delta_{20}s^{\lambda_2}\big(1+\F^\infty_{\ell_2}(\mu_0)\big),\text{ with $\Delta_{20}\!:=\big(L_{22}L_{21}^{-\lambda_2}\big)(1)$}
\end{equation}
and any $0<{\ell}_2<\min(1,\lambda_2^0)$.

Finally, to study the Dulac map $D_3$ of from $\Sigma_1$ to $\Sigma_3$ we make the reflection $\{x_1=y,x_2=x\}$, which brings $-X_\mu$ to \refc{X} with $P_1(x_1,x_2)=-g(x_2,x_1)$ and 
$P_2(x_1,x_2)=-f(x_2,x_1)$. In these coordinates the transverse sections $\Sigma_1$ and $\Sigma_3$ are parametrised by $\sigma_1(s)=(s,1)$ and $\sigma_2(s)=(1,s)$, respectively, and the hyperbolicity ratio of the saddle is $\frac{1}{\lambda_3}=-\big(\frac{P_2}{P_1}\big)(0,0)=-\big(\frac{f}{g}\big)(0,0).$ Hence by \propc{3punts} once again,
\begin{equation}\label{prop1eq3}
 D_3(s)=\Delta_{30}s^{1/\lambda_3}\big(1+\F^\infty_{\ell_3}(\mu_0)\big),\text{ with $\Delta_{30}\!:=\big(L_{32}L_{31}^{-1/\lambda_3}\big)(1)$}
\end{equation}
and any $0<{\ell}_3<\min(1,1/\lambda_3^0)$.
Consequently, from \refc{prop1eq1}, \refc{prop1eq2} and \refc{prop1eq3} we get that
 \begin{align*}
 \mathscr D(s)=\big(D_2\circ D_1-D_3\big)(s)
     &=\Delta_{20}\Delta_{10}^{\lambda_2}s^{\lambda_1\lambda_2}
         \big(1+\F^\infty_{\ell_1}(\mu_0)\big)^{\lambda_2}\big(1+\F^\infty_{\ell_4}(\mu_0)\big)-\Delta_{30}s^{1/\lambda_3}
         \big(1+\F^\infty_{\ell_3}(\mu_0)\big)\\
    &=\Delta_{20}\Delta_{10}^{\lambda_2}s^{\lambda_1\lambda_2}
         \big(1+\F^\infty_{\ell_1}(\mu_0)\big)\big(1+\F^\infty_{\ell_4}(\mu_0)\big)-\Delta_{30}s^{1/\lambda_3}
         \big(1+\F^\infty_{\ell_3}(\mu_0)\big) \\
   &=\Delta_{20}\Delta_{10}^{\lambda_2}s^{\lambda_1\lambda_2}-\Delta_{30}s^{1/\lambda_3}+\F^\infty_{\ell_5}(\mu_0),                  
 \end{align*}
where, by \cite[Lemma A.2]{MV20}, we can take $\ell_4$ and $\ell_5$ to be any numbers such that  $0<\ell_4<\min(\lambda_1^0,\lambda_1^0\lambda_2^0)$ and  $\min(\lambda_1^0\lambda_2^0,1/\lambda_3^0))<\ell_5<\min(\lambda_1^0+\lambda_1^0\lambda_2^0,
1+\lambda_1^0\lambda_2^0,2\lambda_1^0\lambda_2^0,1+1/\lambda_3^0,2/\lambda_3^0)$, respectively. 
Thus, setting $a_1\!:=\Delta_{20}\Delta_{10}^{\lambda_2}$ and $a_2\!:=\Delta_{30}$, we obtain
\begin{equation}\label{prop1eq6}
 \mathscr D(s)=a_1s^{\lambda_1\lambda_2}-a_2s^{1/\lambda_3}+\F^\infty_{\ell_5}(\mu_0)
\end{equation}
and this proves~$(a)$ because each $\Delta_{i0}$ is a strictly positive analytic function by \propc{3punts}. 

Let us proceed next with the proof of the assertions in $(b)$, thus from now on we assume that 
$\lambda_1^0\lambda_2^0\lambda_3^0=1$. Then, setting 
$\alpha=1/\lambda_3-\lambda_1\lambda_2$ and taking any $\ell_6\in (0,\min(1,\lambda_1^0,\lambda_1^0\lambda_2^0))$, from the above expression we get
  \begin{align*}
  \mathscr D(s)
  &=s^{1/\lambda_3}\big(a_1s^{\lambda_1\lambda_2-1/\lambda_3}
       -a_2+\F^\infty_{\ell_6}(\mu_0)\big)\\
  &=s^{1/\lambda_3}\big(a_1(1+\alpha\omega(s;\alpha))-a_2
      +\F^\infty_{\ell_6}(\mu_0)\big)\\
    &=s^{1/\lambda_3}\big(b_1\omega(s;\alpha)+b_2+\F^\infty_{\ell_6}(\mu_0)\big),            
 \end{align*}
where $b_1\!:=\alpha a_1$ and $b_2\!:=a_1-a_2$. Here we also take \defic{defi_comp} into account. Since $a_1$ is an analytic non-vanishing function at $\mu_0$ and $d_1=\lambda_3\alpha$, we obtain the equality $(b_1)=(d_1)$ between ideals over the local ring $\R\{\mu\}_{\mu_0}$.
In order to show that $(b_1,b_2)=(d_1,d_2)$ holds as well we note that, from \refc{prop1eq1}, \refc{prop1eq2} and~\refc{prop1eq3} again,
\begin{align*}
 b_2(\mu)
    &=a_1-a_2=\Delta_{20}\Delta_{10}^{\lambda_2}-\Delta_{30}\\[5pt]
    &=(L_{22}L_{21}^{-\lambda_2}L_{12}^{\lambda_2}L_{11}^{-\lambda_1\lambda_2}
        -L_{32}L_{31}^{-1/\lambda_3})\big|_{u=1}\\[5pt]
    &=L_{32}L_{31}^{-1/\lambda_3}\left(\frac{L_{22}}{L_{32}}\frac{L_{31}^{1/\lambda_3}}{L_{11}^{\lambda_1\lambda_2}}\frac{L_{12}^{\lambda_2}}{L_{21}^{\lambda_2}}-1\right)\Big|_{u=1}
   \\[5pt]
    &=L_{32}L_{31}^{-1/\lambda_3}\left(L_{31}^{\alpha}\left(\frac{L_{22}}{L_{32}}\frac{L_{31}^{\lambda_1\lambda_2}}{L_{11}^{\lambda_1\lambda_2}}\frac{L_{12}^{\lambda_2}}{L_{21}^{\lambda_2}}-1\right)+(L_{31}^\alpha-1)\right)\Big|_{u=1}
         \\[5pt]
    &=\kappa_1(\mu)d_2(\mu)+\alpha(\mu) \kappa_2(\mu),           
\end{align*}
where $\kappa_1=\frac{L_{32}}{L_{31}^{1/\lambda_3}}\big|_{u=1}\frac{e^{d_2}-1}{d_2}$ and $\kappa_2=\frac{L_{32}}{L_{31}^{1/\lambda_3}}\frac{L_{31}^\alpha-1}{\alpha}\big|_{u=1}$ are analytic functions at $\mu_0$ because so it is each~$L_{ij}$ by \cite[Lemma 2.3]{MV21}. 
Thus, on account of $d_1=\lambda_3\alpha$, we get that $b_2=\kappa_1d_2+\kappa_3 d_1$. Hence $(b_1,b_2)=(d_1,d_2)$ over the ring $\R\{\mu\}_{\mu_0}$ since $\kappa_1(\mu_0)> 0$ and we have already shown that $(b_1)=(d_1).$

So far we have proved the first assertion in $(b)$. To show the second one, besides $\lambda_1^0\lambda_2^0\lambda_3^0=1$, we assume $\lambda_1^0>1$, $\lambda_2^0>1$ and $\lambda_3^0<1$. On account of this we can apply point $(3)$ in \propc{3punts} to conclude that
\begin{align}\label{prop1eq61}
 D_1(s)&=s^{\lambda_1}\big(\Delta_{10} +\Delta_{11} s +\F_{\ell_7}^\infty(\mu_0)\big)
 \text{ for any }{\ell_7}\in \big[1,\min(\lambda_1^0,2)\big),\\[3pt]\label{prop1eq7}
 D_2(s)&=s^{\lambda_2}\big(\Delta_{20} +\Delta_{21} s +\F_{\ell_8}^\infty(\mu_0)\big)
 \text{ for any }{\ell_8}\in \big[1,\min(\lambda_2^0,2)\big),
 \intertext{and}\label{prop1eq4}
 D_3(s)&=s^{1/\lambda_3}\big(\Delta_{30} +\Delta_{31} s +\F_{\ell_9}^\infty(\mu_0)\big)
 \text{ for any }{\ell_9}\in \big[1,\min(1/\lambda_3^0,2)\big).
\end{align} 
Here the first order coefficients $\Delta_{10}$, $\Delta_{20}$ and $\Delta_{30}$ are the ones already defined in \refc{prop1eq1}, \refc{prop1eq2} and \refc{prop1eq3}, respectively. With regard to the second order coefficients, only the ones of $D_1$ and $D_3$ are relevant for our purposes, which are given by
\begin{equation}\label{prop1eq5}
 \Delta_{11}\!:=-\frac{\lambda_1\Delta_{10}\hat M_1(1/\lambda_1,1)}{L_{11}(1)}
 \text{ and }
 \Delta_{31}\!:=-\frac{\Delta_{30}\hat M_3(\lambda_3,1)}{\lambda_3 L_{31}(1)}, 
\end{equation}
respectively. In each case, on account of \refc{def_M}, this follows easily from the formula $\Delta_1=\lambda \Delta_0 S_1$ given in \propc{3punts} and taking~$S_1$ in \refc{def_S} particularised to $\sigma_1(s)=(s,1)$. 

From \refc{prop1eq61} and \refc{prop1eq7}, by applying \cite[Lemma A.2]{MV20} we can assert that
\begin{align*}
\big(D_2\circ D_1\big)(s)&=
      s^{\lambda_1\lambda_2}\big(\Delta_{10} +\Delta_{11} s +\F_{\ell_7}^\infty\big)^{\lambda_2}
      \big(\Delta_{20} +\Delta_{21}s^{\lambda_1}\big(\Delta_{10} +\Delta_{11} s +\F_{\ell_7}^\infty\big) 
      +\F_{\ell_{10}}^\infty\big)\\
      &=
      s^{\lambda_1\lambda_2}\big(\Delta_{10}^{\lambda_2} +{\lambda_2}\Delta_{10}^{\lambda_2-1}\Delta_{11} s 
      +\F_{\ell_7}^\infty\big)
      \big(\Delta_{20} +\Delta_{21}s^{\lambda_1}\big(\Delta_{10} +\Delta_{11} s\big) +\F_{\ell_{11}}^\infty 
      +\F_{\ell_{10}}^\infty\big)\\
       &=
      s^{\lambda_1\lambda_2}\big(\Delta_{10}^{\lambda_2} +{\lambda_2}\Delta_{10}^{\lambda_2-1}\Delta_{11} s 
      +\F_{\ell_7}^\infty\big)
      \big(\Delta_{20} +\Delta_{10} \Delta_{21}s^{\lambda_1}
      +\Delta_{11}\Delta_{21}s^{\lambda_1+1}+\F_{\ell_{12}}^\infty\big)
\end{align*} 
for any $\ell_{10}\in [\lambda_1^0,\lambda_1^0\min(\lambda_2^0,2)\big)$ in the first equality, 
any $\ell_{11}\in \big[\lambda_1^0+1,\lambda_1^0+\min(\lambda_1^0,2)\big)$ in the second one
and any $\ell_{12}\in \big[\lambda_1^0+1,\min(2\lambda_1^0,\lambda_1^0+2,\lambda_1^0\lambda_2^0)\big)$ in the third one. Furthermore in the second equality we use that, for any $\eta=\eta(\mu),$
\begin{equation}\label{prop1eq8}
 (1+as+\F_\ell)^\eta=(1+as)^\eta+\F_\ell=1+a\eta s+\F_{2-\varepsilon}+\F_\ell,
\end{equation}
for any $\varepsilon>0,$ which in turn follows noting that 
\begin{align*}
 (1+as+\F_\ell)^\eta-(1+as)^\eta&=\textstyle(1+as)^\eta\big((1+\frac{\F_\ell}{1+as})^\eta-1\big)\\
 &=(1+as)^\eta\big((1+\F_\ell)^\eta-1\big)=(1+as)^\eta\F_\ell=\F_\ell.
\end{align*}
Consequently
\[
 \big(D_2\circ D_1\big)(s)=s^{\lambda_1\lambda_2}
  \big(\Delta_{20}\Delta_{10}^{\lambda_2}+{\lambda_2}\Delta_{11}\Delta_{20}\Delta_{10}^{\lambda_2-1} s
   +\F_{\ell_7}^\infty\big),
\]
where we use again \cite[Lemma A.2]{MV20} taking $\lambda_1^0>1$ and $\lambda_2^0>1$ into account. Hence, from \refc{prop1eq4} and plug in $s^{-\alpha}=1+\alpha\omega(s;\alpha)$ as before, we get
 \begin{align*}
 \mathscr D(s)&=\big(D_2\circ D_1-D_3\big)(s)\\
     &=s^{1/\lambda_3}\Big(
         s^{-\alpha}\big(\Delta_{20}\Delta_{10}^{\lambda_2}+{\lambda_2}\Delta_{11}\Delta_{20}\Delta_{10}^{\lambda_2-1} s
   +\F_{\ell_7}^\infty\big)-\Delta_{30} -\Delta_{31} s -\F_{\ell_9}^\infty\Big)\\
   &=s^{1/\lambda_3}\big((1+\alpha\omega(s;\alpha))\Delta_{20}\Delta_{10}^{\lambda_2}U(s)-\Delta_{30} -\Delta_{31} s -\F_{\ell_9}^\infty\Big),                     
 \end{align*}
where we define
$U(s)=1+{\lambda_2}\Delta_{11}\Delta_{10}^{-1} s+\F_{\ell_7}^\infty$. The application of the formula given in~\refc{prop1eq8} with $\eta=-1$ shows that $U(s)^{-1}=1-\lambda_2\Delta_{11}\Delta_{10}^{-1}s+\F_{\ell_7}$. Thus one can easily verify that the above expression yields to
\[
\mathscr D(s)=s^{1/\lambda_3}U(s)\big(b_1\omega(s;\alpha)+b_2+b_3 s+ \F_{\ell_{13}}^{ \infty}\big)
\]
with $\ell_{13}\in [1,\min(\lambda_1^0,1/\lambda_3^0,2))$ and $b_3\!:=\lambda_2\Delta_{11}\Delta_{30}\Delta_{10}^{-1}-\Delta_{31}$. Let us recall here that $b_1=\alpha a_1=\alpha\Delta_{20}\Delta_{10}^{\lambda_2}$ and $b_2=a_1-a_2=\Delta_{20}\Delta_{10}^{\lambda_2}-\Delta_{30}$, where $a_1$ and $a_2$ are the analytic and strictly positive functions in~\refc{prop1eq6}. On account of the assumptions $\lambda_1^0>1$ and $\lambda_3^0<1$ we have $\lambda_1(\mu)\notin\frac{1}{\N}$ and $\lambda_3(\mu)\notin\N$ for $\mu\approx\mu_0,$ which imply respectively that $\Delta_{11}$ and $\Delta_{31}$ are analytic at $\mu_0$ by \propc{3punts}. Consequently $b_3$ is an analytic function at $\mu_0$. That being said we claim that the equality $(b_1,b_2,b_3)=(d_1,d_2,d_3)$ between ideals over the local ring $\R\{\mu\}_{\mu_0}$ is true. In order to prove this, for the sake of shortness in the next computation we follow the convention that $\kappa$ stands for an analytic function at $\mu_0$ and $\hat\kappa$ stands for an analytic strictly positive function at $\mu_0.$ Some easy computations following this convention yield
\begin{align*}
 b_3&=\Delta_{10}^{-1}({\lambda_2}\Delta_{11}\Delta_{30}-\Delta_{31}\Delta_{10})
    =\Delta_{10}^{-1}({\lambda_2}\Delta_{11}\Delta_{30}-\Delta_{31}\Delta_{10})\\
 &=\textstyle -\Delta_{30}\Big({\lambda_2}\lambda_1\frac{\hat M_1(1/\lambda_1,1)}{L_{11}(1)}
       -\frac{1}{\lambda_3}\frac{\hat M_3(\lambda_3,1)}{L_{31}(1)}\Big)
     =\textstyle -\hat\kappa \Big({\lambda_1}\lambda_2\lambda_3\hat M_1\big(\frac{1}{\lambda_1},1\big)L_{31}(1)
       -\hat M_3(\lambda_3,1)L_{11}(1)\Big)\\
   &=\textstyle -\hat\kappa \Big((1-d_1)\hat M_1\big(\frac{1}{\lambda_1},1\big)L_{31}(1)
       -\hat M_3(\lambda_3,1)L_{11}(1)\Big)=\textstyle\hat\kappa \Big(\hat M_3(\lambda_3,1)L_{11}(1)-\hat M_1\big(\frac{1}{\lambda_1},1\big)L_{31}(1)
       \Big)+\kappa d_1.
\end{align*}
where in the third and fifth equalities we use \refc{prop1eq5} and $d_1\!:=1-{\lambda_1}\lambda_2\lambda_3$, respectively. Hence $b_3=\hat\kappa d_3+\kappa d_1$ since $d_3\!:=\hat M_3(\lambda_3,1)L_{11}(1)-\hat M_1(\frac{1}{\lambda_1},1)L_{31}(1).$ On account of $(b_1)=(d_1)$ and $(b_1,b_2)=(d_1,d_2)$, this shows the validity of the claim and completes the proof of the result.
\end{prova}

\begin{obs}\label{rem_stab}
There are two important observations to be made about \teoc{prop1}:
\begin{enumerate}[$(a)$]
\item
The statement claims that the equalities $(b_1)=(d_1)$, $(b_1,b_2)=(d_1,d_2)$ and $(b_1,b_2,b_3)=(d_1,d_2,d_3)$ between ideals over the local ring $\R\{\mu\}_{\mu_\star}$ are satisfied. As a matter of fact, in the proof we show a stronger property, namely that the following holds:
\[  
 \left(\begin{array}{c} d_1\\ d_2\\ d_3\end{array}\right)
 =\left(\begin{array}{cccc}\kappa_1 & 0 & 0\\ * & \kappa_{2}  & 0\\  * & * & \kappa_{3}\end{array} \right)\left(\begin{array}{c} b_1\\ b_2\\ b_3\end{array}\right),
\]
where all the entries in the matrix are analytic functions on $\Lambda$ and each $\kappa_i$ is strictly positive. 

\item From the dynamical point of view it is interesting to point out that the asymptotic expansion of the displacement map 
$\mathscr D(s;\mu_0)$ given in \teoc{prop1} provides a method to study the stability of the polycycle $\Gamma$. Indeed, taking also the previous observation into account, it shows that
\begin{enumerate}[1.]

\item if $d_1(\mu_0)<0$ (respectively, $>0$) then $\Gamma$ is asymptotically stable (respectively, unstable),

\item if $d_1(\mu_0)=0$ and $d_2(\mu_0)<0$ (respectively, $>0$)
         then $\Gamma$ is asymptotically stable (respectively, unstable), and 
\item if $d_1(\mu_0)=d_2(\mu_0)=0$ and $d_3(\mu_0)<0$ (respectively, $>0$) then $\Gamma$ is asymptotically
         stable (respectively, unstable). 
\end{enumerate}
Of course this is relevant because we have an explicit expression of these functions by \teoc{thmA}. In this regard let us note that the first assertion is well known since $d_1(\mu_0)<0$ is equivalent to require 
that $\lambda_1^0\lambda_2^0\lambda_3^0>1$, while the second assertion was already proved by Gasull \emph{et al.} in \cite{GMM02}, see Theorem~1. On the contrary the third assertion constitutes a new result to the best of our knowledge. 
\end{enumerate}
\end{obs}

We give at this point the precise definition of independence of functions that we use in the present paper.

\begin{defi}\label{indepe}
Consider the functions \map{g_i}{\Lambda}{\R} for $i=1,2,\ldots,k$. The \emph{real variety} $V(g_1,g_2,\ldots,g_k)$ is defined to be the set of $\mu\in\Lambda$ such that $g_i(\mu)=0$ for $i=1,2,\ldots,k.$ We say that $g_1,g_2,\ldots,g_k$ are \emph{independent} at $\mu_{\star}\in V(g_1,g_2,\ldots,g_k)$ if the following conditions are satisfied:\begin{enumerate}[$(1)$]
\item Every neighbourhood of $\mu_{\star}$ contains two points $\mu_1,\mu_2\in V(g_1,\ldots,g_{k-1})$
         such that $g_k(\mu_1)g_k(\mu_2)<0$ (if $k=1$ then we set $V(g_1,\ldots,g_{k-1})=V(0)=\Lambda$ 
         for this to hold).
\item
The varieties $V(g_1,\ldots,g_i)$, $2\leqslant i\leqslant k-1,$ are such that if $\mu_0\in V(g_1,\ldots,g_i)$ 
then every neighbourhood of $\mu_0$ contains two points $\mu_1,\mu_2\in V(g_1,\ldots,g_{i-1})$ such 
        that $g_i(\mu_1)g_i(\mu_2)<0$.
        
\item
If $\mu_0\in V(g_1)$ then every neighbourhood of $\mu_0$ contains two points 
         $\mu_1,\mu_2$ such that $g_{1}(\mu_1)g_{1}(\mu_2)<0.$        
\end{enumerate}
It is clear that if $g_i\in\mathscr C^1(\Lambda)$ for $i=1,2,\ldots,k$ and the gradients $\nabla g_1(\mu_{\star}),\nabla g_2(\mu_{\star})\dots,\nabla g_k(\mu_{\star})$ are linearly independent vectors of $\R^{N+1}$ 
then there exists a neighbourhood $U_\star$ of $\mu_\star$ such that the restrictions of 
 $g_1,g_2,\ldots,g_k$ to $U_\star$ are independent at $\mu_{\star}$.
 \end{defi}

\begin{lem}\label{ideals}
Suppose that the equalities $(c_1,\ldots,c_k)=(d_1,\ldots,d_k)$ between ideals over the local ring $\R\{\mu\}_{\mu_\star}$ hold for $k=1,2,\ldots,n$, where $\mu_\star\in V(c_1,\ldots,c_n)=V(d_1,\ldots,d_n).$ Then $c_1,\ldots,c_n$ are independent at $\mu_\star$ if, and only if, $d_1,\ldots,d_n$ are independent at $\mu_\star$.
\end{lem}

\begin{prova}
Let us assume for instance that $c_1,\ldots,c_n$ are independent at $\mu_\star$ and prove that then $d_1,\ldots,d_n$ are also independent. To this aim we note that the equalities $(c_1,\ldots,c_k)=(d_1,\ldots,d_k)$ for $k=1,2,\ldots,n$ imply the existence of two triangular matrices 
$A=(a_{ij})$ and $B=(b_{ij})$ with coefficients in $\R\{\mu\}_{\mu_\star}$ such that 
\[
 \left(\begin{array}{c} d_1\\ \vdots\\ d_n\end{array}\right)
 =\left(\begin{array}{cccc}a_{11} & 0 & \cdots & 0\\ * & a_{22} & \cdots & 0\\ \vdots & \vdots & \ddots & \vdots\\ * &\cdots & * & a_{nn}\end{array} \right)\left(\begin{array}{c} c_1\\ \vdots\\ c_n\end{array}\right)
\text{ and }
\left(\begin{array}{c} c_1\\ \vdots\\ c_n\end{array}\right)=\left(\begin{array}{cccc} b_{11} & 0 & \cdots & 0\\ * & b_{22} & \cdots & 0\\ \vdots & \vdots & \ddots & \vdots\\ * &\cdots & * & b_{nn}\end{array} \right)\left(\begin{array}{c} d_1\\ \vdots\\ d_n\end{array}\right).
\]
Clearly $R=(r_{ij})\!:=BA$ is also a triangular matrix with coefficients in the local ring $\R\{\mu\}_{\mu_\star}$ and
\[
\left(\begin{array}{c} c_1\\ \vdots\\ c_n\end{array}\right)
=\left(\begin{array}{cccc} r_{11} & 0 & \cdots & 0\\ * & r_{22} & \cdots & 0\\ \vdots & \vdots & \ddots & \vdots\\ * &\cdots & * & r_{nn}\end{array} \right)\left(\begin{array}{c} c_1\\ \vdots\\ c_n\end{array}\right).
\]
We claim that, since $c_1,\ldots,c_n$ are independent at $\mu_\star$, then $r_{kk}(\mu_{\star})=1$ for all $k=1,2,\ldots,n.$ The fact that this is true for $k=1$ follows easily by continuity. Let us prove by contradiction that this is also true for $k\geqslant 2.$ So assume that $r_{kk}(\mu_\star)\neq 1$ for some $k\in \{2,\ldots,n\}.$ Then the equality $c_k=r_{k1}c_1+\ldots+r_{kk}c_k$ implies that $c_k=\alpha_1c_1+\ldots+\alpha_{k-1}c_{k-1}$ where each $\alpha_{i}\!:=\frac{r_{ki}}{1-r_{kk}}$ is an analytic function at $\mu_\star.$ This clearly contradicts the assumption that $c_1,\ldots,c_n$ are independent at $\mu_\star$ (see \defic{indepe}). Hence the claim is true and, consequently, $\text{det}(R)=\text{det}(A)\text{det}(B)=1$ at $\mu=\mu_\star.$ This shows, in particular, that $A$ is an invertible matrix in the local ring $\R\{\mu\}_{\mu_\star}$ and so there exists a neighbourhood $U$ of $\mu_\star$ such that $a_{kk}(\mu)\neq 0$ for all $\mu\in U$ and $k=1,2,\ldots,n$. On account of this, the fact that $d_1,\ldots,d_n$ are independent at~$\mu_\star$ follows easily noting that if we take any two points $\mu_1,\mu_2\in U\cap V(c_1,\ldots,c_{i-1})=U\cap V(d_1,\ldots,d_{i-1})$ verifying $c_i(\mu_1)c_i(\mu_2)<0$ then we have that $d_i(\mu_1)d_i(\mu_2)=a_{ii}(\mu_1)a_{ii}(\mu_2)c_i(\mu_1)c_i(\mu_2)<0$. This completes the proof of the result.
\end{prova}

\begin{prooftext}{Proof of \teoc{thmA}.}
Let us fix any $\mu_0\in\Lambda$ and set $\lambda_i^0\!:=\lambda_i(\mu_0)$ for $i=1,2,3.$ Recall that the limit cycles of $X_\mu$ near $\Gamma$ are in one to one correspondence with the isolated positive zeros of 
 \[
  \mathscr D(s;\mu)=\big(D_2\circ D_1-D_3\big)(s;\mu)
 \]
near $s=0.$ If $d_1(\mu_0)=1-\lambda_1^0\lambda_2^0\lambda_3^0$ is not zero then by applying $(a)$ in \teoc{prop1} we have that, for any $\ell_1\in\big(\min(\lambda_1^0\lambda_2^0,1/\lambda_3^0),\min(\lambda_1^0+\lambda_1^0\lambda_2^0,
1+\lambda_1^0\lambda_2^0,2\lambda_1^0\lambda_2^0,1+1/\lambda_3^0,2/\lambda_3^0)\big)$,
\begin{equation*}
 \mathscr D(s;\mu)=a_1(\mu)s^{\lambda_1\lambda_2}+a_2(\mu)s^{1/\lambda_3}+\F_{\ell_1}^\infty(\mu_0),
\end{equation*}
where $a_1$ and $a_2$ are analytic and strictly positive functions on $\Lambda.$ Thus
\begin{align*}
 \lim_{(s,\mu)\to (0,\mu_0)}s^{-\lambda_1\lambda_2}\mathscr D(s;\mu)=a_1(\mu_0)
 \text{ in case that $\lambda_1^0\lambda_2^0<\lambda_3^0$}\\
 \intertext{and}
\lim_{(s,\mu)\to (0,\mu_0)}s^{-1/\lambda_3}\mathscr D(s;\mu)=a_2(\mu_0) 
 \text{ in case that $\lambda_1^0\lambda_2^0>\lambda_3^0$}.
\end{align*}
Since $a_i(\mu_0)\neq 0$ for $i=1,2$, this implies the existence of an open neighbourhood $U$ of $\mu_0$ and $\varepsilon>0$ small enough such that $\mathscr D(s;\mu)\neq 0$ for all $\mu\in U$ and $s\in (0,\varepsilon)$ when $\lambda_1^0\lambda_2^0\lambda_3^0\neq 1.$ Hence $\mathrm{Cycl}\big((\Gamma,X_{\mu_0}),X_\mu\big)=0$ and the assertion in $(a)$ is true.

In order to prove $(b)$ we note that if $d_1(\mu_0)=1-\lambda_1^0\lambda_2^0\lambda_3^0$ is equal to zero then, by $(b)$ in \teoc{prop1}, we can assert that, for any 
$\ell_2\in(0,\min(1,\lambda_1^0,\lambda_1^0\lambda_2^0))$,
\begin{equation}\label{thAeq2}
 s^{-1/\lambda_3}\mathscr D(s;\mu)=b_1(\mu)\,\omega\big(s;\alpha(\mu)\big)+b_2(\mu)+\F_{\ell_2}^\infty(\mu_0),
\end{equation}
where $\alpha=1/\lambda_3-\lambda_1\lambda_2$ and $b_1$ and $b_2$ are analytic functions at $\mu_0$ such that 
$(b_1)=(d_1)$ and $(b_1,b_2)=(d_1,d_2)$ over the ring $\R\{\mu\}_{\mu_0}$. The assumptions in this case imply that $\mathscr D(s;\mu_0)\not\equiv 0$, $b_1(\mu_0)=0$ and, thanks to \lemc{ideals}, that $b_1$ is independent at $\mu_0.$ 
Thus, given any $\varepsilon>0$, there exists $s_1\in (0,\varepsilon)$ such that $\mathscr D(s_1;\mu_0)\neq 0.$ Let us assume for instance that $\mathscr D(s_1;\mu_0)> 0,$ the other case follows verbatim. Then, thanks to $(1)$ in \defic{indepe}, there exists $\mu_1\approx\mu_0$ with $b_1(\mu_1)<0$ and, by continuity, such that $\mathscr D(s_1;\mu_1)> 0.$ Moreover by applying Lemmas A.3 and A.4 in \cite{MV19},
\begin{equation}\label{thAeq1}
 Z_1(s;\mu)\!:=\frac{s^{-1/\lambda_3}\mathscr D(s;\mu)}{\omega\big(s;\alpha(\mu)\big)}=b_1(\mu)+\frac{b_2(\mu)}{\omega\big(s;\alpha(\mu)\big)}+\F_{\ell_2-\delta}^\infty(\mu_0)\to \kappa(\mu)\text{ as $s\to 0$,}
\end{equation}
where
\[
 \kappa(\mu)\!:=b_1(\mu)-b_2(\mu)\min(\alpha(\mu),0)).
\]
Here we use that $1/{\omega(s;\alpha(\mu))}\in\F_{-\delta}^\infty(\mu_0)$ for any $\delta>0$ and that $\lim_{s\to 0}1/{\omega(s;\alpha)}=\max(-\alpha,0)$ by assertions~$(a)$ and $(b)$ in \cite[Lemma A.4]{MV19}, respectively. Note on the other hand that, by $(b)$ in \teoc{prop1}, $b_1=\alpha a_1$ and $b_2=a_1-a_2$, where each $a_i$ is an analytic strictly positive function. Thus $\kappa(\mu)=\alpha(\mu)a_2(\mu)$ if $\alpha(\mu)<0$ and $\kappa(\mu)=\alpha(\mu)a_1(\mu)$ if $\alpha(\mu)\geqslant 0.$
Therefore, since $b_1=\alpha a_1$, we can write
\begin{equation}\label{thAeq3}
 \kappa=b_1\eta\text{ with $\eta>0$.} 
\end{equation}
Hence, on account of $b_1(\mu_1)<0$, we can assert that $\kappa(\mu_1)<0,$ which in turn, from \refc{thAeq1}, guarantees the existence of some $s_2\in (0,s_1)$ such that $Z_1(s_2;\mu_1)<0.$ Thus $\mathscr D(s_1;\mu_1)\mathscr D(s_2;\mu_1)<0$ and, by continuity, $\mathscr D(\hat s;\mu_1)=0$ for some $\hat s\in (s_2,s_1)\subset (0,\varepsilon).$ This shows that $\mathrm{Cycl}\big((\Gamma,X_{\mu_0}),X_\mu\big)\geqslant 1$ as desired.

Let us prove $(c)$ next, i.e., that if $d_2(\mu_0)\neq 0$ then $\mathrm{Cycl}\big((\Gamma,X_{\mu_0}),X_\mu\big)\leqslant 1$. Note first that to this end we can also assume that $d_1(\mu_0)=0$, otherwise $\mathrm{Cycl}\big((\Gamma,X_{\mu_0}),X_\mu\big)=0.$ Consequently the expression in \refc{thAeq2} is valid. That being said, the proof will 
follow by applying the well-known derivation-division algorithm. In doing so, from \refc{thAeq1} and the fact that $\partial_s\omega(s;\alpha)=-s^{-\alpha-1},$
\begin{equation*}
 \partial_sZ_1(s;\mu)
 =\frac{b_2(\mu)}{s^{\alpha(\mu)+1}\omega^2\big(s;\alpha(\mu)\big)}+\F_{\ell_2-1-\delta}^\infty(\mu_0),
\end{equation*}
where the flatness of the remainder follows from $(f)$ in Lemma A.3 of \cite{MV19}. Therefore
\begin{equation*}
 s^{\alpha(\mu)+1}\omega^2\big(s;\alpha(\mu)\big)\partial_sZ_1(s;\mu)=b_2(\mu)+\F_{\ell_2-4\delta}^\infty(\mu_0)\to b_2(\mu_0)
 \text{ as $(s,\mu)\to (0,\mu_0).$} 
\end{equation*}
Recall at this point that $(b_1)=(d_1)$ and $(b_1,b_2)=(d_1,d_2)$ over the local ring $\R\{\mu\}_{\mu_0}$ by \teoc{prop1}. Thus, the assumptions $d_1(\mu_0)=0$ and $d_2(\mu_0)\neq 0$ imply that $b_2(\mu_0)\neq 0$. Accordingly, on account of the above limit and by Bolzano's Theorem, we obtain $\varepsilon>0$ such that if $\|\mu-\mu_0\|<\varepsilon$ then $Z_1(\,\cdot\,;\mu)$, and so $\mathscr D(\,\cdot\,;\mu)$, has at most one zero for $s\in (0,\varepsilon)$, multiplicities taking into account. Hence $\mathrm{Cycl}\big((\Gamma,X_{\mu_0}),X_\mu\big)\leqslant 1$ and~$(c)$ follows. 

Let us turn next to the proof of $(d)$, in which the assumptions are $\mathscr R(\,\cdot\,;\mu_0)\not\equiv \text{Id}$ and that $d_1$ and $d_2$ vanish and are independent at $\mu_0$. Then $\mathscr D(\,\cdot\,;\mu_0)\not\equiv 0$ and, due to $(b_1)=(d_1)$ and $(b_1,b_2)=(d_1,d_2)$ once again, $b_1$ and $b_2$ vanish and are independent at $\mu_0$ by \lemc{ideals}. Thus, given any $\varepsilon>0$, there exists $s_1\in (0,\varepsilon)$ such that, for instance, $\mathscr D(s_1;\mu_0)<0.$ Then, by continuity and condition $(1)$ in \defic{indepe}, there exists $\mu_1\approx\mu_0$ such that $b_2(\mu_1)>0$, $b_1(\mu_1)=0$ and $\mathscr D(s_1;\mu_1)<0.$ Hence, from \refc{thAeq2}, 
\[
s^{-1/\lambda_3}\mathscr D(s;\mu_1)=b_2(\mu_1)+\F_{\ell_2}^\infty(\mu_0)\to b_2(\mu_1)\text{ as $s\to 0,$}
\]
which shows the existence of $s_2\in (0,s_1)$ such that $\mathscr D(s_2;\mu_1)>0.$ For the same reasons we can choose $\mu_2\approx\mu_1$ satisfying $\mathscr D(s_1;\mu_2)<0$ and $\mathscr D(s_2;\mu_2)>0$ together with $b_1(\mu_2)<0.$ 
Then, from~\refc{thAeq1} and~\refc{thAeq3}, $\lim_{s\to 0}Z_1(s;\mu_2)=b_1(\mu_2)\eta(\mu_2)<0$ and so there exists $s_3\in (0,s_2)$ verifying that $\mathscr D(s_3;\mu_2)<0.$ By continuity there exist $\hat s_1,\hat s_2\in (0,\varepsilon)$ with $D(\hat s_1;\mu_2)=D(\hat s_2;\mu_2)=0.$ Accordingly $\mathrm{Cycl}\big((\Gamma,X_{\mu_0}),X_\mu\big)\geqslant 2$.

From now on, in order to prove $(e)$ and $(f)$, we assume $\lambda_1^0<1$, $\lambda_2^0>1$ and $\lambda_3^0>1$. Then by applying \teoc{prop1}, for any $\ell_3\in \big(1,\min(2,\lambda_1^0,1/\lambda_3^0)\big)$,
\begin{equation}\label{thAeq4}
 \mathscr D(s;\mu)=\big(b_1(\mu)\,\omega\big(s;\alpha(\mu)\big)+b_2(\mu)+b_3(\mu)s+\F_{\ell_3}^\infty(\mu_0)\big)s^{1/\lambda_3}U(s;\mu),
\end{equation}
where $b_3$ is an analytic function at $\mu_0$ verifying that $(b_1,b_2,b_3)=(d_1,d_2,d_3)$ 
over the ring $\R\{\mu\}_{\mu_0}$ and $U$ is an analytic function such that $U(0;\mu_0)=1$. Hence
\begin{equation}\label{thAeq5}
 Z_2(s;\mu)\!:=\frac{\mathscr D(s;\mu)}{s^{1/\lambda_3}\omega(s;\alpha(\mu))U(s)}=b_1(\mu)+b_2(\mu)\frac{1}{\omega(s;\alpha(\mu))}+b_3(\mu)\frac{s}{\omega(s;\alpha(\mu))}+\F_{\ell_3-\delta}^\infty(\mu_0),
\end{equation}
where we use once again that $1/{\omega(s;\alpha(\mu))}\in\F_{-\delta}^\infty(\mu_0)$ for any $\delta>0$. Note furthermore that, for $\mu\approx\mu_0$, the positive zeros of $\mathscr D(\,\cdot\,;\mu)$ and $Z_2(\,\cdot\,;\mu)$ near $s=0$ are in one to one correspondence because $\frac{1}{s^{1/\lambda_3}\omega(s;\alpha(\mu))}$ tends to $+\infty$ as $(s,\mu)\to (0,\mu_0).$
That being stablished we begin first with the proof of assertion~$(e)$ and to this aim, besides $d_3(\mu_0)\neq 0,$ we can also suppose $d_1(\mu_0)=d_2(\mu_0)=0$, otherwise $\mathrm{Cycl}\big((\Gamma,X_{\mu_0}),X_\mu\big)\leqslant 1$ by~$(a)$ or~$(b)$, which already have been proved. In this case, since $(b_1)=(d_1)$, $(b_1,b_2)=(d_1,d_2)$ and  $(b_1,b_2,b_3)=(d_1,d_2,d_3)$, it turns out that $b_3(\mu_0)\neq 0.$ As in the proof of $(c)$, we will apply to steps of the derivation-division algorithm in \refc{thAeq5}. In doing so we obtain that
\[
 Z_3(s;\mu)\!:=s^{\alpha(\mu)+1}\omega^2(s;\alpha(\mu))\partial_sZ_2(s;\mu)=
 b_2(\mu)+b_3(\mu)\big(s+s^{\alpha(\mu)+1}\omega(s;\alpha(\mu))\big)+\F_{\ell_3-4\delta}^\infty(\mu_0),
\]
where the flatness of the remainder follows by applying Lemmas A.3 and A.4 in \cite{MV19} as before and we use that $\partial_s\omega(s;\alpha)=-s^{-\alpha-1}.$ Note also that the positive zeros of $Z_3(\,\cdot\,;\mu)$ and $\partial_sZ_2(\,\cdot\,;\mu)$ near $s=0$ are in one to one correspondence for $\mu\approx\mu_0$ because $\omega(s;\alpha(\mu))$ tends to $+\infty$ as $(s,\mu)\to (0,\mu_0).$ Finally
\[
 \frac{\partial_sZ_3(s;\mu)}{s^{\alpha(\mu)}\omega(s;\alpha(\mu))}=(\alpha(\mu)+1)b_3(\mu)+\F_{\ell_3-1-6\delta}^\infty(\mu_0)\to b_3(\mu_0)\neq 0\text{ as $(s,\mu)\to (0,\mu_0)$.}
\]
By applying twice Bolzano's Theorem, we can assert the existence of some $\varepsilon>0$ such that if $\|\mu-\mu_0\|<\varepsilon$ then $Z_2(\,\cdot\,;\mu)$, and so $\mathscr D(\,\cdot\,;\mu)$, has at most two zeros for $s\in (0,\varepsilon)$, multiplicities taking into account. Hence $\mathrm{Cycl}\big((\Gamma,X_{\mu_0}),X_\mu\big)\leqslant 2$, which proves $(e)$. 

Finally, in order to prove $(f)$ we suppose that $\mathscr R(\,\cdot\,;\mu_0)\not\equiv \text{Id}$ and that $d_1$, $d_2$ and $d_3$ vanish and are independent at $\mu_0$. Consequently $\mathscr D(\,\cdot\,;\mu_0)\not\equiv 0$ and, due to the equality between the corresponding ideals over the local ring, $b_1$, $b_2$ and $b_3$ vanish and are independent at~$\mu_0$ by \lemc{ideals}. Thus, given any $\varepsilon>0$, there exists $s_1\in (0,\varepsilon)$ such that, for instance, $\mathscr D(s_1;\mu_0)<0.$ Then, by continuity and condition $(1)$ in \defic{indepe}, there exists $\mu_1\approx\mu_0$ such that $b_3(\mu_1)>0$, $b_1(\mu_1)=b_2(\mu_1)=0$ and $\mathscr D(s_1;\mu_1)<0.$ Hence, from \refc{thAeq4}, 
\[
\frac{\mathscr D(s;\mu_1)}{s^{1+1/\lambda_3}U(s;\mu)}=b_3(\mu_1)+\F_{\ell_3-1}^\infty(\mu_0)\to b_3(\mu_1)>0\text{ as $s\to 0,$}
\]
which shows the existence of $s_2\in (0,s_1)$ such that $\mathscr D(s_2;\mu_1)>0.$ For the same reasons we can choose $\mu_2\approx\mu_1$ satisfying $\mathscr D(s_1;\mu_2)<0$ and $\mathscr D(s_2;\mu_2)>0$ together with $b_1(\mu_2)=0$ and $b_2(\mu_2)<0.$ Accordingly, from~\refc{thAeq4} again, 
\[
\frac{\mathscr D(s;\mu_2)}{s^{1/\lambda_3}U(s;\mu)}=b_2(\mu_2)+b_3(\mu_2)s+\F_{\ell_3-1}^\infty(\mu_0)\to b_2(\mu_2)<0\text{ as $s\to 0,$}
\]
which shows the existence of $s_3\in (0,s_2)$ such that $\mathscr D(s_3;\mu_2)<0.$ In the final step we take 
$\mu_3\approx\mu_2$ satisfying $\mathscr D(s_1;\mu_3)<0$ and $\mathscr D(s_2;\mu_3)>0$ and $\mathscr D(s_3;\mu_3)<0$ together with $b_1(\mu_3)>0.$ 
Then, from~\refc{thAeq1} and~\refc{thAeq3}, $\lim_{s\to 0}s^{-1/\lambda_3}Z_1(s;\mu_3)=b_1(\mu_3)\eta(\mu_3)>0$ and so there exists $s_4\in (0,s_3)$ such that $\mathscr D(s_3;\mu_3)>0.$ By continuity there exist $\hat s_1,\hat s_2,\hat s_3\in (0,\varepsilon)$ with $D(\hat s_i;\mu_3)=0$ for $i=1,2,3.$ Hence $\mathrm{Cycl}\big((\Gamma,X_{\mu_0}),X_\mu\big)\geqslant 3$ and this completes the proof of the result.
\end{prooftext}

\section{Applications}\label{appli}

We begin this section by revisiting in \teoc{ex1} a family of Kolmogorov differential systems that was first studied in \cite{GMM02}, where the authors (following the notation in our statement) prove that if $\mu_0=(a_0,p_0,q_0)$ verifies $p_0+q_0=0$ and $a_0\neq 0$ then $\mathrm{Cycl}\big((\Gamma,X_{\mu_0}),X_\mu\big)\geqslant 1$, cf. assertion $(b)$. 

\begin{theo}\label{ex1}
Consider the family of Kolmogorov differential systems 
\[
 X_\mu\quad\sist{x(1+x+x^2+axy+py^2),}{y(-1-y+qx^2+axy-y^2),}
\]
where $\mu=(a,p,q)\in\R^3$ with $p<-1$ and $q>1$ and let us fix any $\mu_0=(a_0,p_0,q_0).$ Then, compactifying $X_\mu$ to the Poincaré disc, the boundary of the first quadrant is a polycycle $\Gamma$ such that:
\begin{enumerate}[$(a)$]
\item $\mathrm{Cycl}\big((\Gamma,X_{\mu_0}),X_\mu\big)=0$ if $p_0+q_0\neq 0.$

\item $\mathrm{Cycl}\big((\Gamma,X_{\mu_0}),X_\mu\big)=1$ if $p_0+q_0=0$ 
         and $a_0\neq 0$.
         
\item The return map of $X_{\mu_0}$ along $\Gamma$ is the identity 
         if, and only if, $a_0=p_0+q_0=0$. In this case $\Gamma$ is the
         outer boundary of the period annulus of a center at $(x_0,y_0)$ with 
         $x_0=y_0=-\frac{1+\sqrt{-3-4p_0}}{2(1+p_0)}$ that foliates the first quadrant and, moreover, 
         $\mathrm{Cycl}\big((\Gamma,X_{\mu_0}),X_\mu\big)=1$.   
 
\end{enumerate}

\noindent On the other hand the vector field $X_{\mu}$ has a unique singularity $Q_{\mu}=(\upsilon_1,\upsilon_2)$ in the first quadrant, which is either a focus or a center, and has trace equal to $\tau(\mu)=\upsilon_1+2\upsilon_1^{2}+2a\upsilon_1\upsilon_2-\upsilon_2-2\upsilon_2^{2}$. Furthermore the following holds:
\begin{enumerate}[$(a)$]\setcounter{enumi}{3}            
         
\item If $\tau(\mu_0)\neq 0$ then $\mathrm{Cycl}\big((Q_{\mu_0},X_{\mu_0}),X_\mu\big)=0$ and a sufficient
         condition for $\tau(\mu_0)\neq 0$ to hold is that $p_0+q_0=0$ and $a_0\neq 0.$
         
\item If $\tau(\mu_0)=0$ and $p_0+q_0\neq 0$ then $Q_{\mu_0}$ is a weak focus of order 1 and 
         $\mathrm{Cycl}\big((Q_{\mu_0},X_{\mu_0}),X_\mu\big)=1.$   
         
\item If $\tau(\mu_0)= 0$ and $p_0+q_0= 0$ then $a_0=0$. In addition $Q_{\mu_0}$ is a center if, and only if,
         $p_0+q_0=a_0=0$ and, in this case, $\mathrm{Cycl}\big((Q_{\mu_0},X_{\mu_0}),X_\mu\big)=1.$        

\end{enumerate}
Finally it is not possible a simultaneous bifurcation of limit cycles from $\Gamma$ and $Q_\mu$.

\end{theo}

\begin{prova}
The assertions in $(a)$ and $(b)$ follow directly by applying \teoc{thmA}. Indeed, 
in this case, following the notation in \refc{kol}, $f(x,y)=1+x+x^2+axy+py^2$ and $g(x,y)=-1-y+qx^2+axy-y^2$, so that
\[
 f_2(x,y)=x^2+axy+py^2\text{ and }g_2(x,y)=qx^2+axy-y^2.
\]
Taking this into account, together with $p<-1$ and $q>1,$ one can easily check that the assumptions {\bf H1} and {\bf H2} are verified. As a matter of fact the first assumption holds not only for $z>0$ but for all $z\in \R$, and this implies that the boundary of each quadrant is a monodromic polycycle for the compactified vector field. Hence, by the Poincaré-Bendixson theorem (see \cite{Roussarie} for instance), there exists at least one singularity of~$X_{\mu}$ inside each one of the four quadrants. Due to $\deg(f)=\deg(g)=2$, by Bézout's theorem there exists exactly one in each quadrant. From now on we denote the singularity of $X_\mu$ in the first quadrant by $Q_{\mu}$. That being said, the hyperbolicity ratios of the saddles at $\Gamma$ are $\lambda_1=\frac{1}{q-1}$, $\lambda_2=-(p+1)$ and $\lambda_3=1.$ Consequently the first assertion follows from $(a)$ in \teoc{thmA} because 
\begin{equation}\label{ex1eq1}
d_1(\mu)=1-\lambda_1\lambda_2\lambda_3=\frac{p+q}{q-1}.
\end{equation}
The second assertion will follow by applying $(b)$ and $(c)$ in \teoc{thmA}. To show this we first recall that
\[
 d_2(\mu)=\lambda_2\log\left(\frac{L_{12}}{L_{21}}\right)\!(1)+\log\left(\frac{L_{22}}{L_{32}}\right)\!(1)+\lambda_1\lambda_2\log\left(\frac{L_{31}}{L_{11}}\right)\!(1)
\]
and this leads us to the computation of the following improper integrals:
\begin{align*}
 \Delta_1(\mu)\!:=\log\left(\frac{L_{12}}{L_{21}}\right)\!(1)&=\int_0^1\left(\left(\frac{f_2}{f_2-g_2}\right)(1,z)+\lambda_1-\left(\frac{g_2}{g_2-f_2}\right)(z,1)-\frac{1}{\lambda_2}\right)\frac{dz}{z}\\[5pt]
 \Delta_2(\mu)\!:=\log\left(\frac{L_{22}}{L_{32}}\right)\!(1)&=\int_0^1\left(\left(\frac{g-f}{g}\right)(0,1/z)+\lambda_2
 -\left(\frac{f}{g}\right)(0,z)-\frac{1}{\lambda_3}\right)\frac{dz}{z}\\[5pt]
  \Delta_3(\mu)\!:=\log\left(\frac{L_{31}}{L_{11}}\right)\!(1)&=\int_0^1\left(\left(\frac{g}{f}\right)(z,0)+\lambda_3
  -\left(\frac{f-g}{f}\right)(1/z,0)-\frac{1}{\lambda_1}\right)\frac{dz}{z}
\end{align*}
These expressions have to be computed assuming that $p+q=0,$ i.e., $\lambda_1\lambda_2=1.$ In doing so we obtain that
\begin{align*}
 &\Delta_1(a,p,-p)=\frac{2a}{1-q}\int_0^1\frac{zdz}{z^2+1}=\frac{a\pi}{2(p+1)}\\
 \intertext{and}
 &\Delta_2(a,p,-p)=-\Delta_3(a,p,-p)=-(p+1)\int_0^1\frac{zdz}{z^2+z+1}=-\frac{(p+1)\pi}{3\sqrt{3}}.
\end{align*}
Therefore $d_2(a,p,-p)=-\frac{a\pi}{2}$ is zero if, and only if, $a=0.$ Taking this into account, the combination of $(b)$ and $(c)$ in \teoc{thmA} shows that $\mathrm{Cycl}\big((\Gamma,X_{\mu_0}),X_\mu\big)=1$ for any $\mu_0=(a_0,p_0,-p_0)$ with $a_0\neq 0$, as desired. It is important to remark for the forthcoming analysis that by applying the Weierstrass Division Theorem (see for instance \cite{Greuel,Krantz}) we can assert that
\begin{equation}\label{ex1eq2}
 d_2(\mu)=-\frac{a\pi}{2}+(p+q)h(\mu)
\end{equation}
for some analytic function $h.$

Next we proceed with the proof of $(c)$.To this aim we fix any $\mu_0=(a_0,p_0,q_0)$ and apply \teoc{prop1}, which gives the asymptotic expansion of $\mathscr D(s;\mu)$ at $s=0$ for $\mu\approx\mu_0$. This result, taking~\refc{ex1eq1} and~\refc{ex1eq2} into account, shows that if $\mathscr D(s;\mu_0)\equiv 0$ then $a_0=p_0+q_0=0$. In order to prove the converse observe that if $\mu_0=(0,p_0,-p_0)$ then the vector field $X_{\mu_0}$ writes as
\[
 \sist{x(1+x+x^2+p_0y^2),}{-y(1+y+p_0x^2+y^2).}
\]
One can easily check that $Q_{\mu_0}$, the only singularity of $X_{\mu_0}$ in the first quadrant, is a weak focus at the point $(x_0,y_0)$ with $x_0=y_0=-\frac{1+\sqrt{-3-4p_0}}{2(1+p_0)}$. Furthermore, setting $\sigma(x,y)=(y,x)$, it turns out that $\sigma^\star X_{\mu_0}=-X_{\mu_0}$ and so the vector field is reversible with respect to the straight line $y=x.$ Hence $Q_{\mu_0}$ is a center and a straightforward application of the Poincaré-Bendixson theorem shows that its period annulus fills the first quadrant, which in particular implies that $\mathscr D(s;\mu_0)\equiv 0$. 

So far we have proved that the return map of $X_{\mu_0}$ along $\Gamma$ is the identity if, and only if, $\mu_0=(0,p_0,-p_0).$ Our next task is to show that, in this case, $\mathrm{Cycl}\big((\Gamma,X_{\mu_0}),X_\mu\big)=1$. With this aim in view we apply $(b)$ in \teoc{prop1}, which shows that if $\mu\approx\mu_0$ then 
\begin{equation}\label{ex1eq7}
 \mathscr D(s;\mu)=\big(b_1(\mu)\,\omega\big(s;\alpha(\mu)\big)+b_2(\mu)+r(s;\mu)\big)s^{1/\lambda_3},
\end{equation}
where $\alpha=1/\lambda_3-\lambda_1\lambda_2$, $r\in \F_{\ell}^\infty(\mu_0)$ with $\ell\in \big(0,\min(1,\frac{-1}{p+1})\big)$ and, in addition,
\[
 (b_1)=(d_1)\text{ and }(b_1,b_2)=(d_1,d_2)
\]
over the local ring $\R\{\mu\}_{\mu_0}$. Consequently if $\mu=(a,p,q)$ satisfies $a=p+q=0$ then  $b_1(\mu)=b_2(\mu)=0$ and $r(s;\mu)\equiv 0.$ Furthermore, since the vectors $\nabla d_1(\mu_0)$ and $\nabla d_2(\mu_0)$ are linearly independent, see~\refc{ex1eq1} and~\refc{ex1eq2}, the above equalities between ideals show that this is also the case of $\nabla b_1(\mu_0)$ and $\nabla b_2(\mu_0).$ We can thus take $(\eta_1,\eta_2,\eta_3)\notin \langle\nabla b_1(\mu_0),\nabla b_2(\mu_0)\rangle$ and define $b_3(\mu)\!:=\eta_1a+\eta_2(p-p_0)+\eta_3(q+p_0)$ so that $\nu=\Psi(\mu)\!:=\big(b_1(\mu),b_2(\mu),b_3(\mu)\big)$ is a local analytical change of coordinates in a neighbourhood
 of $\mu=\mu_0$. Note that $\Psi$ maps $\mu_0$ to $0_3\!:=(0,0,0)$ and $\{a=p+q=0\}$ to $\{\nu_1=\nu_2=0\}$ and in addition 
\begin{equation*}
 \mathscr R_1(s;\nu)\!:=\left.s^{-1/\lambda_3}\mathscr D(s;\mu)\right|_{\mu=\Psi^{-1}(\nu)}=\nu_1\,\omega(s;\hat\alpha)+\nu_2+\hat r(s;\nu),
\end{equation*}
where $\hat\alpha=\hat\alpha(\nu)\!:=\alpha \big(\Psi^{-1}(\nu)\big)$ and 
$\hat r(s;\nu)\!:=r \big(s;\Psi^{-1}(\nu)\big)\in\F_{\ell}^\infty(0_3).$ The key point is that $\hat r(s;0,0,\nu_3)\equiv 0$ implies, thanks to \cite[Lemma 4.1]{MV22}, that $\hat r(s;\nu)=\nu_1h_1(s;\nu)+\nu_2h_2(s;\nu)$ with $h_i\in\F_{\ell}^\infty(0_3).$
Accordingly
\begin{equation}\label{ex1eq10}
 \mathscr R_1(s;\nu)=\nu_1\big(\omega(s;\hat\alpha)+h_1(s;\nu)\big)+\nu_2\big(1+h_2(s;\nu)\big).
\end{equation}
Observe that if $s\to 0^+$ and $\nu\to (0,0,0)$ then the factor multiplying $\nu_1$ tends to $+\infty$, whereas the factor multiplying $\nu_2$ tends to $1.$ Here we use \defic{defi2} and that $\lim_{(s,\alpha)\to (0,0)}\omega(s;\alpha)=+\infty.$ We claim that there exists $s_0>0$ and an open neighbourhood~$U$ of $\nu=(0,0,0)$ such that 
\[
 \mathscr R_2(s;\nu)\!:=\frac{\mathscr R_1(s;\nu)}{\omega(s;\hat\alpha)+h_1(s;\nu)}=\nu_1+\nu_2\frac{1+h_2(s;\nu)}{\omega(s;\hat\alpha)+h_1(s;\nu)}
\]
has at most one zero on $(0,s_0)$, counted with multiplicities, for all $\nu=(\nu_1,\nu_2,\nu_3)\in U$ with $\nu_1^2+\nu_2^2\neq 0.$ This will imply that $\mathrm{Cycl}\big((\Gamma,X_{\mu_0}),X_\mu\big)\leqslant 1$ because 
$\mathscr R_2(s;0,0,\nu_3)\equiv 0$, so that it has not any isolated zero. The claim is clear in case that $\nu
_2=0.$ To tackle the case $\nu_2\neq 0$ we compute the derivative with respect to~$s$ to obtain that
\begin{align*}
 \mathscr R_2'(s;\nu)&=\nu_2\partial_s\left(
 \frac{1+\F_{\ell}^\infty}{\omega(s;\hat\alpha)+\F_{\ell}^\infty} \right)
 =\nu_2\partial_s\left(\frac{1+\F_{\ell}^\infty}{\omega(s;\hat\alpha)(1+\F_{\ell-\varepsilon}^\infty)}
 \right)\\[6pt]
 &=\nu_2\partial_s\left(\frac{1+\F_{\ell-\varepsilon}^\infty}{\omega(s;\hat\alpha)}\right)
 =\frac{\nu_2}{s^{\hat\alpha+1}\omega^2(s;\hat\alpha)}(1+\F_{\ell-\varepsilon}^\infty)
  +\frac{\nu_2}{\omega(s;\hat\alpha)}\F_{\ell-\varepsilon-1}^\infty\\[6pt]
  &=\frac{\nu_2}{s^{\hat\alpha+1}\omega^2(s;\hat\alpha)}\big(1+\F_{\ell-\varepsilon}^\infty+s^{\hat\alpha+1}\omega(s;\hat\alpha)\F_{\ell-\varepsilon-1}^\infty\big)
  =\frac{\nu_2}{s^{\hat\alpha+1}\omega^2(s;\hat\alpha)}\big(1+\F_{\ell-3\varepsilon}^\infty\big).
\end{align*}
Here, in the second equality we apply first assertion $(c)$ of Lemma A.4 in~\cite{MV19}
 to get that $1/\omega(s;\hat\alpha)\in\F^\infty_{-\varepsilon}(0_3)$ for all $\varepsilon>0$ small enough, due to $\hat\alpha(0_3)=0$, and use next that $\F_{-\varepsilon}^\infty\F_{\ell}^\infty\subset\F_{\ell-\varepsilon}^\infty$ from $(g)$ of Lemma~A.3 in~\cite{MV19}
. In the third equality,
on account of $\frac{1}{1+s}-1\in\F_1^\infty$ and by $(h)$ of Lemma~A.3 in~\cite{MV19}
, we use first the inclusion $\frac{1}{1+\F^\infty_{\ell-\varepsilon}}\subset 1+\F^\infty_{\ell-\varepsilon}$. Then, by using $(d)$ and $(g)$ of Lemma~A.3 in~\cite{MV19}
, we expand the numerator to get that $(1+\F_{\ell}^\infty)(1+\F_{\ell-\varepsilon}^\infty)\subset 1+\F_{\ell-\varepsilon}^\infty$. Next, in the fourth equality we use that $\partial_s\omega(s;\alpha)=s^{-\alpha-1}$ and assertion $(f)$ of Lemma A.3 in \cite{MV19}
 to deduce that $\partial_s\F_{\ell-\varepsilon}^\infty\subset\F_{\ell-\varepsilon-1}^\infty$. Finally in the last equality we apply $(c)$ of Lemma~A.4 in \cite{MV19}
 to get that $s^{\hat\alpha+1}\omega(s;\hat\alpha)\in\F^\infty_{1-2\varepsilon}$ and we use again that $\F_{1-2\varepsilon}^\infty\F_{\ell-\varepsilon-1}^\infty\subset\F_{\ell-3\varepsilon}^\infty$.
On account of \defic{defi2} we can assert the existence of some $s_0\in (0,1)$ and a neighbourhood $U$ of~$\nu=(0,0,0)$ such that $\mathscr R_2'(s;\nu
)\neq 0$ for all $s\in (0,s_0)$ and $\nu
\in U$ with $\nu_2\neq 0.$ Hence the application of Rolle's theorem shows that the claim is true for $\nu_2\neq 0$ as well. So far we have proved that $\mathrm{Cycl}\big((\Gamma,X_{\mu_0}),X_\mu\big)\leqslant 1$. The fact that this upper bound is attained follows by applying the assertion in $(b)$ taking $\mu_0=(a_0,p_0,-p_0)$ with $a_0\approx 0$ but different from zero. This completes the proof of $(c)$. 

Let us turn now to the proof of the assertions regarding the singularity of $X_\mu$ at $Q_\mu.$ The approach here is rather standard and the technical difficulty is that we do not dispose of a feasible expression of the coordinates of $Q_\mu$. To overcome this problem we shall parametrise the family of vector fields more conveniently.  For reader's convenience we summarise the chain of reparametrisations that we shall perform:
\[
 \mu=(a,p,q)\to (a,p,\varepsilon)\to (\upsilon_1,\upsilon_2,\varepsilon)\to (\upsilon_1,\upsilon_2,\tau)\to \hat\mu=(\ka,\kb,\tau).
\]
For the first one we simply introduce $\varepsilon=p+q.$ In the second one we take the coordinates of the singular point $Q_\mu=(\upsilon_1,\upsilon_2)$ as new parameters, i.e., we isolate $a$ and $p$ from  
\begin{equation}\label{ex1eq4}
  \left\{\!
   \begin{array}{l}
    1+\upsilon_1+\upsilon_1^2+a\upsilon_1\upsilon_2+p\upsilon_2^2=0, \\[2pt] 
    1+\upsilon_2+(p-\varepsilon)\upsilon_1^2-a\upsilon_1\upsilon_2+\upsilon_2^2=0,
   \end{array}
  \right.
\end{equation}
to obtain
\[
a=-\frac{\upsilon_1^{2} \upsilon_2^{2}\varepsilon+\upsilon_1^{4}-\upsilon_2^{4}+\upsilon_1^{3}-\upsilon_2^{3}+\upsilon_1^{2}-\upsilon_2^{2}}{\upsilon_1  \upsilon_2  \left(\upsilon_1^{2}+\upsilon_2^{2}\right)}
\text{ and }
p=\frac{\upsilon_1^{2} \varepsilon -\upsilon_1^{2}-\upsilon_2^{2}-\upsilon_1 -\upsilon_2 -2}{\upsilon_1^{2}+\upsilon_2^{2}}.
\]
In this respect we point out that $\upsilon_1$ and $\upsilon_2$ are strictly positive because $Q_\mu$ is inside the first quadrant for all admissible $\mu$.
More important, the map $\varphi:(a,p,\varepsilon)\mapsto(\upsilon_1,\upsilon_2,\varepsilon)$ is smooth and, taking \refc{ex1eq4} into account, injective. The smoothness follows by the Inverse Function Theorem since one can check that the determinant of the Jacobian of $(\upsilon_1,\upsilon_2,\varepsilon)\mapsto(a,p,\varepsilon)$ is non-zero 
at the image by $\varphi$ of any admissible parameter.
Then one can check that the trace of the Jacobian of the vector field at $(x,y)=(\upsilon_1,\upsilon_2)$ is   
\begin{align*}
 xf_x(x,y)+yg_y(x,y)\big|_{(x,y)=(\upsilon_1,\upsilon_2)}&=
 \upsilon_1+2\upsilon_1^{2}+2a\upsilon_1\upsilon_2-\upsilon_2-2\upsilon_2^{2}\Big|_{a=a(\upsilon_1,\upsilon_2,\varepsilon)}\\[3pt]
  &=
 -{\frac {2\varepsilon\upsilon_1^2\upsilon_2^2+ ( \upsilon_1-\upsilon_2) ( \upsilon_1+2+\upsilon_{{2}} ) ( \upsilon_1+\upsilon_2)}{\upsilon_2^{2}+\upsilon_1^{2}}},
\end{align*}
and we introduce $\tau=\tau(\upsilon_1,\upsilon_2,\varepsilon)$ isolating $\varepsilon$ from 
\begin{equation}\label{ex1eq6}
 2\varepsilon\upsilon_1^2\upsilon_2^2+ ( \upsilon_1-\upsilon_2) ( \upsilon_1+2+\upsilon_{{2}} ) ( \upsilon_1+\upsilon_2)=\tau.
\end{equation}
In other words, $\tau$ is (up to a non-vanishing factor) the trace of the vector field. Finally, for convenience, we define 
$\ka=\frac{\upsilon_1-\upsilon_2}{2}$ and $\kb=\frac{\upsilon_1+\upsilon_2}{2}.$ Observe then that $\{p+q=a=0\}$ becomes $\{\ka=\tau=0\}.$ In what follows, setting $\hat\mu=(\ka,\kb,\tau)$ for shortness, we denote the vector field by $X_{\hat\mu}$. Let us also remark that the map $\mu\mapsto\hat\mu$ is smooth and injective as a consequence of the previous discussion.

At this point we claim that $Q_\mu$ is either a focus or a center. To show this we will check that the discriminant $D_\mu$ of the characteristic polynomial of the Jacobian matrix of $X_\mu$ at $Q_\mu$ is strictly negative for all admissible parameter. Indeed, one can verify that $D_\mu$ expressed in terms of $(\upsilon_1,\upsilon_2,\varepsilon)$ can be written as
 \[
  D_\mu=\frac{4(\upsilon_1\upsilon_2)^4\varepsilon^2+A_1(\upsilon_1,\upsilon_2)\varepsilon+A_0(\upsilon_1,\upsilon_2)}{(\upsilon_1^2+\upsilon_2^2)^2},
 \]
where $A_i$ are polynomials of degree 7. Thus $D_\mu=0$ gives two roots $\varepsilon=\hat\varepsilon_i(\upsilon_1,\upsilon_2)$, for $i=1,2,$ that one can check to be well-defined continuous functions on $V\!:=\{(\upsilon_1,\upsilon_2)\in\R^2:\upsilon_1>0,\upsilon_2>0\}.$ To see the claim we first prove that, for $i=1,2,$
\begin{equation}\label{ex1eq5}
 (p+1)(q-1)\big|_{\varepsilon=\hat\varepsilon_i}>0\text{ for all $(\upsilon_1,\upsilon_2)\in V.$}
\end{equation}
This implies that $D_\mu$ can not vanish at an admissible parameter due to the assumptions $p<-1$ and $q>1.$ In this regard we note that the product $(p+1)(q-1)$ expressed in terms of $(\upsilon_1,\upsilon_2,\varepsilon)$ is given by
\[
 (p+1)(q-1)=\frac{(\upsilon_1\upsilon_2)^2\varepsilon^2+B_1(\upsilon_1,\upsilon_2)\varepsilon+B_0(\upsilon_1,\upsilon_2)}{(\upsilon_1^2+\upsilon_2^2)^2},
\] 
where $\deg(B_1)=3$ and $\deg(B_0)=2.$ A computation shows that the \emph{resultant} (see \cite[Chapter 3]{Cox} for instance) between the numerators of~$D_\mu$ and $(p+1)(q-1)$ with respect to $\varepsilon$ is a polynomial in~$\upsilon_1$ and~$\upsilon_2$ with all the coefficients being natural numbers. Consequently the resultant does not vanish on $V$ and, accordingly, $(p+1)(q-1)|_{\varepsilon=\hat\varepsilon_i}\neq 0$ for all $(\upsilon_1,\upsilon_2)\in V.$ Thus, since $V$ is arc-connected and the function $(\upsilon_1,\upsilon_2)\mapsto (p+1)(q-1)|_{\varepsilon=\hat\varepsilon_i}$ is continuous on $V$, it suffices to verify~\refc{ex1eq5} at some particular choice of parameter. For instance, taking $\upsilon_1=\upsilon_2=1$ we obtain that $(p+1)(q-1)|_{\varepsilon=\hat\varepsilon_i}=92$ for $i=1,2.$ 
Therefore $D_\mu\neq 0$ at any admissible parameter. Thus, exactly as before, since $\mu\mapsto D_\mu$ is continuous and the set of admissible parameters is arc-connected, the claim will follow once we verify its validity at some particular parameter. For instance the choice $\mu=(0,-2,2)$ yields $D_\mu=-\frac{67+25\sqrt{5}}{2}<0.$

We proceed now with the study of the cyclicity of $Q_\mu$. The fact that $\mathrm{Cycl}\big((Q_{\mu_0},X_{\mu_0}),X_\mu\big)=0$ when $\tau(\mu_0)\neq 0$ is well-known. On the other hand, if $p+q=0$ then $\varepsilon=0$ and 
\[ 
 \tau(\upsilon_1,\upsilon_2,0)= ( \upsilon_1-\upsilon_2) ( \upsilon_1+2+\upsilon_{{2}} ) ( \upsilon_1+\upsilon_2),
\] 
which vanishes at an admissible parameter if, and only if, $\upsilon_1=\upsilon_2$. For this to happen, see \refc{ex1eq4} with $\varepsilon=0$, it is necessary that $a=0.$ This proves the validity of $(d)$. 

We shall next solve the center-focus problem in the family. With this aim in view, taking a local transversal section at $Q_\mu$ we consider the displacement map $\mathscr D(s;\hat\mu)$, which extends analytically to $s=0$, so that we can compute its Taylor's expansion 
\begin{equation}\label{ex1eq8}
 \mathscr D(s;\hat\mu)=\eta_1(\hat\mu)s+\eta_2(\hat\mu)s^2+\eta_3(\hat\mu)s^3+sR(s;\hat\mu),
\end{equation}
where the remainder $R$ is $\op(s^2)$. Recall that the trace of $X_{\hat\mu}$ at $Q_{\hat\mu}$ is equal to $\tau u_1(\hat\mu)$, where $u_1$ is a unity. The coefficients $\eta_i$ are called the \emph{Lyapunov quantities} of the focus. We have in particular (see for instance \cite[p. 94]{RS}) that $\eta_1(\hat\mu)=e^{\tau u_1(\hat\mu)}-1=\tau u_2(\hat\mu)$, where $u_2$ is again a unity. Since the first nonzero coefficient of the expansion is the coefficient of an odd power of $s$, see \cite[p. 94]{RS} again, we get that $\eta_2(\hat\mu)=\tau \ell_1(\hat\mu)$ for some analytic function~$\ell_1$. In order to obtain $\eta_3$ we shall appeal to the well-known relation between the Lyapunov and focus quantities which, following the notation in \cite[Theorem 6.2.3]{RS}, we denote by~$g_{ii}$. The first ones are the coefficients in 
the Taylor's expansion of the displacement map that we already introduced, while the second ones are the obstructions for the existence of a first integral. It occurs that $\eta_{2i+1}-\pi g_{ii}\in (g_{11},\ldots,g_{i-1,i-1})$ and, more important for our purposes, that $\eta_3=\pi g_{11}.$ On account of this we can compute $g_{11}$ instead of $\eta_3$, which is easier to obtain, and in doing so (see \cite[p. 29]{Christopher06}) we get that
\[
 \eta_3(\hat\mu)\big|_{\tau=0}=\pi{\frac { 2\ka\left( \ka-\kb \right)  \left( \kb+1
 \right)  \left( 4+34\kb+29{\kb}^{2} +8{\kb}^{3}-3{\ka}^{2}\right) }{ 3\left( \ka+\kb \right) ^{3} \left( \ka+2
\kb+2{\kb}^{2}+2{\ka}^{2}+2 \right) ^{2}}}.
\] 
In this respect we claim that $\eta_3(\hat\mu)\big|_{\tau=0}=\ka h(\ka,\kb)$ with $h(\ka,\kb)\neq 0$ in case that $|\ka|<\kb$, which corresponds to the admissible values $\upsilon_1,\upsilon_2>0$ due to $\ka=\frac{\upsilon_1-\upsilon_2}{2}$ and $\kb=\frac{\upsilon_1+\upsilon_2}{2}.$ Indeed, it is clear that the factor $(\ka-\kb)(\kb+1)$ does not vanish inside the admissible set, while the other one does not vanish neither because
\[
  4+34\kb+29{\kb}^{2}+8{\kb}^{3}-3{\ka}^{2}>4+34\kb+26{\kb}^{2}+8{\kb}^{3}>0,
\] 
where the first inequality follows using that $|\ka|<\kb$ and the second one the fact that $\kb>0.$  Hence the claim is true. Therefore, if $Q_{\mu}$ is a center then $\tau=\ka=0$, and the assertion in $(c)$ shows that these two conditions are also sufficient because $\{p+q=a=0\}=\{\tau=\ka=0\}.$ Observe moreover that we can write $\eta_3(\hat\mu)=\tau \ell_2(\hat\mu)+\ka h(\ka,\kb)$ for some analytic function $\ell_2.$ On the other hand, due to $R(s;\hat\mu)|_{\ka=\tau=0}\equiv 0$, we can also write $R=\tau R_{1}+\ka R_{2}$ with $R_i\in\op(s^2)$ and, accordingly,
\begin{equation}\label{ex1eq3}
 \mathscr D(s;\mu)=\tau s\big(u_2+\ell_1s+\ell_2s^2+R_{1})+\ka s\big(hs^2+R_{2}\big).
\end{equation}
Note that if $\tau(\mu)=0$ and $\varepsilon=p+q\neq 0$ then $\ka=\frac{\upsilon_1-\upsilon_2}{2}$ must be different from zero because otherwise, from~\refc{ex1eq6}, we would get that $\varepsilon=0.$ Consequently, due to $h(\ka,\kb)\neq 0$ for all admissible $\ka$ and $\kb,$ the equality in \refc{ex1eq3} implies $\mathrm{Cycl}\big((Q_{\mu_0},X_{\mu_0}),X_\mu\big)\leqslant 1$ in case that $\tau(\mu_0)=0$ and $p_0+q_0\neq 0.$ The fact that this upper bound is attained follows by means of an easy perturbative argument using that $\partial_\varepsilon\tau(\mu)=-2\upsilon_1^2\upsilon_2^2\neq 0.$ This proves the validity of the assertion in $(e).$

In order to prove $(f)$ note that if $\tau(\mu)=0$ and $\varepsilon=p+q=0$ then, from \refc{ex1eq6}, 
$\ka=\frac{\upsilon_1-\upsilon_2}{2}=0.$ Hence, from \refc{ex1eq4}, $2a\upsilon_1\upsilon_2=0$, which implies $a=0$ and shows the first assertion. That being stablished, we have already proved that $Q_{\mu_0}$ is a center if, and only if $p_0+q_0=a_0=0.$ We show next that, in this case, $\mathrm{Cycl}\big((Q_{\mu_0},X_{\mu_0}),X_\mu\big)\leqslant 1$. Indeed, 
since $u_2$ is a unity we can consider
\[
  \mathscr D_1(s;\hat\mu)\!:=\frac{\mathscr D(s;\hat\mu)}{s(u_2+\ell_1s+\ell_2s^2+R_{1})}=\tau+
  \ka\frac{hs^2+R_{2}}{u_2+\ell_1s+\ell_2s^2+R_{1}}.
\] 
The upper bound for the cyclicity of $Q_{\mu_0}$ in the center case will follow once we prove that there exist $s_0>0$ and an open neighbourhood $U$ of $(\ka,\kb,\tau)=(0,\hat\kb,0)$ such that 
$\mathscr D_1(s;\hat\mu)$ has at most one zero on $(0,s_0)$, counted with multiplicities, for all $\hat\mu\in U$ with $\ka^2+\tau^2\neq 0.$ Recall in this regard that $\mathscr D_1(s;\hat\mu)|_{\ka=\tau=0}\equiv 0$ and that, on account of $(c)$, $\hat\kb=-\frac{1+\sqrt{-3-4p_0}}{2(1+p_0)}>0$ is the first component of $Q_{\mu_0}.$ The idea to show this is exactly the same as in the proof of $(c)$ but with less technicalities because the involved functions are analytic at $s=0$. The desired property is evident when $\ka=0.$ In case that $\ka\neq 0$ we compute the derivative of $\mathscr D_1$ with respect to $s$ to obtain 
\[
 \partial_s\mathscr D_1(s;\hat\mu)=\ka s \big(2h/u_2+\op(1)\big).
\]
Since $h(\ka,\kb)\neq 0$ in case that $|\ka|<\kb$ and $u_2=u_2(\hat\mu)$ is a unity, the existence of the desired $s_0>0$ and the open neighbourhood $U$ follows by Rolle's theorem. So far we have proved that $\mathrm{Cycl}\big((Q_{\mu_0},X_{\mu_0}),X_\mu\big)\leqslant 1$ if $p_0+q_0=a_0=0$. The fact that this upper bound is attained follows noting that we can take $\mu=(a,p,q)$ with $\tau(\mu)= 0$ and $p+q\neq 0$ arbitrarily close to $\mu_0=(0,p_0,-p_0)$ and apply then the assertion in $(e).$ This proves $(f)$. 

Let us turn now to the proof of the last assertion in the statement. Observe in this respect that the combination of $(a)$ and $(b)$ together with $(d)$ and $(e)$ shows that a simultaneous bifurcation of limit cycles from $\Gamma$ and $Q_\mu$ can only occur if we perturb some $\mu_\star=(a_\star,p_\star,q_\star)$ with $a_\star=p_\star+q_\star=0.$ We shall prove by contradiction that this is neither possible. So assume that for each $n\in\N$ there exist $\mu_n=(a_n,p_n,q_n)$ and two limit cycles $\gamma_n$ and $\gamma_n'$ of the vector field $X_{\mu_n}$ in the first quadrant such that the Hausdorff distances $d_H(\gamma_n,\Gamma)$ and $d_H(\gamma_n',Q_{\mu_n})$
tend to zero and $\mu_n$ tends to $\mu_\star$ as $n\to +\infty.$ Let us consider the asymptotic expansion of the displacement map of $X_\mu$ at the polycycle $\Gamma$ that we compute in \refc{ex1eq7} and denote it by $\mathscr D_p(s;\mu)$. We also consider its Taylor's expansion near the focus $Q_\mu$ given in \refc{ex1eq3} and denote it by $\mathscr D_c(s';\mu).$ Then the assumption implies the existence of two sequences $s_n\to 0^+$ and $s_n'\to 0^+$ such that $\mathscr D_p(s_n;\mu_n)=0$ and $\mathscr D_c(s'_n;\mu_n)=0$ for all $n\in\N.$ We claim that the first equality implies that 
\begin{equation}\label{ex1eq8}
 \lim_{n\to +\infty}\frac{p_n+q_n}{a_n}=0.
\end{equation}
Indeed, from \refc{ex1eq10} we have that 
\[
 \mathscr R_1(s_n;\nu)\big|_{\nu=\Psi(\mu_n)}=b_1(\mu_n)\,\big(\omega(s_n;\alpha(\mu_n))
 +h_1(s_n;\nu)\big)+b_2(\mu_n)\big(1+h_2(s;\nu)\big)\big|_{\nu=\Psi(\mu_n)}=0\text{ for all $n\in\N.$}
\]
Thus, due to $\lim_{n\to +\infty}\omega (s_n;\alpha(\mu_n))=+\infty$ and $h_i\in\F_{\ell}^\infty(0_3)$, we obtain that 
$\frac{b_2(\mu_n)}{b_1(\mu_n)}\to -\infty$ as $n\to +\infty$.
Moreover, since $(b_1)=(d_1)$ and $(b_1,b_2)=(d_1,d_2)$ with $\nabla d_1$ and $\nabla d_2$ independent at $\mu=\mu_\star,$ we can write
\[
 \frac{b_2}{b_1}=\frac{\kappa_1d_1+\hat\kappa_2d_2}{\hat\kappa_3d_1}
\]
with $\hat\kappa_i(\mu_\star)\neq 0$ and, consequently, $\lim_{n\to +\infty}\frac{d_2(\mu_n)}{d_1(\mu_n)}=\infty$. This, on account of \refc{ex1eq1} and \refc{ex1eq2}, gives the limit in~\refc{ex1eq8} and so the claim is true. Recall on the other hand that in order to study the displacement map near the focus $Q_\mu$ we use a more convenient parametrisation given by $\hat\mu\!:=(\ka,\kb,\tau)=\phi(\mu)$. That being said, setting $(\ka_n,\kb_n,\tau_n)\!:=\phi(a_n,p_n,q_n)$, similarly as we argue before, the fact that $\mathscr D_c(s'_n;\mu_n)=0$ for all $n\in\N$ implies from \refc{ex1eq3} that 
\begin{equation}\label{ex1eq9}
 \lim_{n\to +\infty}\frac{\tau_n}{\ka_n}=0.
\end{equation}
Let us remark that here we also take into account that $u_2$ is a unity. We next arrive to contradiction showing that \refc{ex1eq8} and \refc{ex1eq9} cannot hold simultaneously. Indeed, one can verify that, setting $\kb_\star=-\frac{1+\sqrt{-3-4p_\star}}{2(1+p_\star)},$
\begin{align*}
   \left.\frac{p_n+q_n}{a_n}\right|_{\mu_n=\phi^{-1}(\hat\mu_n)}&=4\,\frac{\ka_n^2+\kb_n^2}{\ka_n^2-\kb_n^2}\,
  \frac{\tau_n+\ka_n\kb_n(\kb_n+2)}{2\tau_n-\ka_n(2\kb_n+1)(\ka_n^2+\kb_n^2)}\\[6pt]
  &=4\,\frac{\ka_n^2+\kb_n^2}{\ka_n^2-\kb_n^2}
  \,\frac{\tau_n/\ka_n+\kb_n(\kb_n+2)}{2\tau_n/\ka_n-(2\kb_n+1)(\ka_n^2+\kb_n^2)}
  \to\frac{4(\kb_\star+2)}{\kb_\star(2\kb_\star+1)}\neq 0\text{ as $n\to+\infty$.}
\end{align*}
Here, in addition to \refc{ex1eq9}, we use that if $p+q$ and $a$ tend to zero then $\ka\to 0$ and $\kb\to\kb_\star,$ where $\kb_\star$ is precisely the first component of the center at $Q_{\mu_\star}$ (which is in the diagonal of the first quadrant). This shows that \refc{ex1eq8} and \refc{ex1eq9} cannot occur simultaneously, which yields to the desired contradiction and finishes the proof of the result.
\end{prova}

The following is our second example of application of \teoc{thmA}. In this case the family of Kolmogorov's systems is five-parametric and to the best of our knowledge it has not been studied previously. 

\begin{theo}\label{ex2}
Consider the family of Kolmogorov differential systems 
\[
 X_\mu\quad\sist{x(c+x^2+axy-(p+1)y^2),}{y(-1+(q+1)x^2+(a-\bb)xy-y^2),}
\]
where $\mu=(a,\bb,c,p,q)\in\R^5$ with $c>0$, $p>0$, $q>0$ and $\bb<2\sqrt{pq}$ and let us fix any $\mu_0=(a_0,\bb_0,c_0,p_0,q_0).$ 
Then there exists a unique singular point $Q_\mu$ in the first quadrant, which is either a center, a focus or a node. Moreover, compactifying~$X_\mu$ to the Poincaré disc, the boundary of the first quadrant is a polycycle $\Gamma$ such that:

\begin{enumerate}[$(a)$]
\item $\mathrm{Cycl}\big((\Gamma,X_{\mu_0}),X_\mu\big)=0$ if $p_0-c_0q_0\neq 0.$

\item $\mathrm{Cycl}\big((\Gamma,X_{\mu_0}),X_\mu\big)=1$ if $p_0-c_0q_0=0$ 
         and $2c_0q_0a_0-(c_0q_0-c_0+1)\bb_0\neq 0$.

\item The return map of $X_{\mu_0}$ along $\Gamma$ is the identity 
         if, and only if, $p_0-c_0q_0=2c_0q_0a_0-(c_0q_0-c_0+1)\bb_0=0$. 
         In this case $Q_{\mu_0}$ is a center with first integral 
         \[
           H(x,y)=\frac{q_0(x^2+c_0(y^2+1))-\bb_0 xy}{(xy^{c_0})^\frac{2}{c_0q_0+c_0+1}},
         \]
         which foliates the first quadrant. Moreover $\Gamma$ is the outer boundary of its period annulus and, 
         in addition, 
         $\mathrm{Cycl}\big((\Gamma,X_{\mu_0}),X_\mu\big)=1$.   
 
\end{enumerate}
\end{theo}

\begin{obs}
In contrast to the family of Kolmogorov's  cubic systems studied in \teoc{ex1},
for the family in \teoc{ex2} there exist parameters $\mu_0$ with $d_1(\mu_0)=0$ and $d_2(\mu_0)\neq 0$, so that 
$\mathrm{Cycl}\big((\Gamma,X_{\mu_0}),X_\mu\big)=1$, and satisfying additionally that the unique singular point $Q_{\mu_0}$ in the first quadrant is a non-degenerate node. Hence, for appropiate $\mu\approx\mu_0$ we will have a limit cycle $\gamma_\mu$ with a non-monodromic singular point~$Q_\mu$ as unique singularity in its interior. For instance, the choice $\mu_0=(-800.01,-900.99999,1000,1,0.001)$ leads to this phenomenon with $Q_{\mu_0}=(0.1,10)$. A similar occurrence is observed in~\cite[p. 203]{ADL} to take place in the family of cubic Lienard systems studied in \cite{DumLi}. 
\end{obs}

\begin{prooftext}{Proof of \teoc{ex2}.}
In this case, following the notation in \refc{kol}, we have that
\[
 f(x,y;\mu)\!:=c+x^2+axy-(p+1)y^2
 \text{ and }
 g(x,y;\mu)\!:=-1+(q+1)x^2+(a-\bb)xy-y^2.
\] 
Since $f(z,0;\mu)=c+z^2$, $g(0,z;\mu)=-1-z^2$ and $(f_2-g_2)(1,z;\mu)=-q+\bb z-pz^2$,
one can check that the hypothesis $\mathbf{H1}$ and $\mathbf{H2}$ are satisfied for the admissible parameters, i.e., $c>0$,
$p>0$, $q>0$ and $\bb<2\sqrt{pq}$. Moreover the hyperbolicity ratios are 
\begin{equation}\label{ex2eq0}
 \text{$\lambda_1=1/q$, $\lambda_2=p$ and 
$\lambda_3=1/c$.}
\end{equation}
Then $\Gamma$ is a polycycle and by applying the Poincar\'e-Bendixson theorem we deduce the existence of at least one singular point of $X_\mu$ in the first quadrant. We claim that there exists exactly one. In order to show this we suppose that $(\upsilon_1,\upsilon_2)$ is a singular point of~$X_\mu$ in the first quadrant and solve $f(\upsilon_1,\upsilon_2;\mu)=0$ and $g(\upsilon_1,\upsilon_2;\mu)=0$ for $a$ and $\bb$ as a function of $c,$ $p,$ $q,$ $\upsilon_1$ and $\upsilon_2$. In doing so we obtain that 
\[
 a={\frac {p\upsilon_2^{2}-\upsilon_1^{2}+\upsilon_{{2}}^{2}-c}{\upsilon_
{{1}}\upsilon_{{2}}}}\text{ and }\bb={\frac {p \upsilon_{{2}}^{2}+q\upsilon_{{1}}^{2}-c-
1}{\upsilon_{{1}}\upsilon_{{2}}}}.
\]
The substitution of these values in $f+cg$, which is homogeneous of degree 2 in $x$ and $y$, yields
\[
 \left.\big(f+cg\big)(x,y;\mu)\right|_{y=rx}=
 \frac{x^2(\upsilon_2-r\upsilon_1)\big(r\upsilon_2(c+p+1)+\upsilon_1(cq+c+1)\big)}{\upsilon_1\upsilon_2}.
\]
It is clear then that the vanishing of the above numerator provides the possible values of $r$ such that $X_\mu$ has a singular point at the straight-line $y=rx,$ namely,
\[
 r_1=\frac{\upsilon_2}{\upsilon_1}>0\text{ and }r_2=-\frac{\upsilon_1(cq+c+1)}{\upsilon_2(c+p+1)}<0.
\]
Since $x\mapsto f(x,r_1x)=c+x^2(1+ar_1-(p+1)r_1^2)$ vanishes at $x=\upsilon_1>0$, it must have another real zero, which has to be negative. Therefore $X_\mu$ has exactly one singular point in the first quadrant and exactly one singular point in the third quadrant, showing in particular the validity of the claim. An easy computation shows that the determinant of the Jacobian of $X_\mu$ at $Q_\mu=(\upsilon_1,\upsilon_2)$ is equal to $2\upsilon_1^2(cq+c+1)+2\upsilon^2_2(p+c+1)>0$, so that it can be a center, a focus or a node.

So far we have proved that the first assertion in the statement is true. Let us turn to the proof of the assertions in $(a)$, $(b)$ and $(c)$. The first one follows from 
$(a)$ in \teoc{thmA} because 
\begin{equation}\label{ex2eq1}
 d_1(\mu)=1-\lambda_1\lambda_2\lambda_3=\frac{cq-p}{cq}.
\end{equation}
The second assertion will follow by applying $(b)$ and $(c)$ in \teoc{thmA}. 
In this regard let us recall that
\begin{equation}\label{ex2d2}
d_2(\mu)= \lambda_2\log\left(\frac{L_{12}}{L_{21}}\right)\!(1)+\log\left(\frac{L_{22}}{L_{32}}\right)\!(1)+\lambda_1\lambda_2\log\left(\frac{L_{31}}{L_{11}}\right)\!(1).
\end{equation}
On account of the definition of each $L_{ij}$, see \refc{def_M}, we easily obtain that
\begin{align}\notag
\log L_{11}(1)&=\frac{1+c(q+1)}{2c}\log(c+1),\\\label{ex2eq6}
\log L_{31}(1)&=\frac{1+c(q+1)}{2c}\log\left(1/c+1\right),\\
\log L_{22}(1)&=\log L_{32}(1)=\frac{\log 2}{2}(p+c+1).\notag
\end{align}
Moreover
\[
\log L_{12}(1)=\frac{1}{q}\int_0^1\frac{mz+n}{-pz^2+\bb z-q}dz\text{ and }
\log L_{21}(1)=\frac{1}{p}\int_0^1\frac{mz+n'}{-qz^2+\bb z-p}dz,
\]
where
\begin{equation}\label{ex2eq3}
 m\!:=-(pq+p+q),\text{ }n\!:=qa+\bb\text{ and }n'\!:=p(b-a)+b.
\end{equation}
The explicit integration of these functions leads to several cases depending on the parameters. To avoid this we note that
\[
 \frac{mz+n}{-pz^2+\bb z-q}=-\frac{m}{2p}\frac{-2pz+b}{-pz^2+\bb z-q}+\frac{1}{2p}\frac{mb+2np}{-pz^2+\bb z-q}, 
\]
so that
\begin{equation*}
\log L_{12}(1)=-\frac{m}{2pq}\log\left(\frac{p+q-\bb}{q}\right)+\frac{m\bb+2np}{2pq}\int_0^1\frac{dz}{-pz^2+bz-q}.
\end{equation*}
It is clear that the same formula holds for $\log L_{21}(1)$ replacing $p,$ $q$ and $n$ by $q$, $p$ and $n'$, respectively.
On account of this and the fact that, from \refc{ex2eq3}, $mb+2n'q=-mb-2np$, we get
\begin{equation}\label{ex2eq4}
 \log\left(\frac{L_{12}}{L_{21}}\right)\!(1)=\log L_{12}(1)-\log L_{21}(1)=
-\frac{m}{2pq}\log\left(\frac{p}{q}\right)+(m\bb+2np)\Phi(\mu),
\end{equation}
where
\[
\Phi(\mu)\!:=\frac{1}{2pq}\int_0^1\left(\frac{1}{-pz^2+\bb z-q}+\frac{1}{-qz^2+\bb z-p}\right)dz.
\]
Notice, and this is the key point in the forthcoming arguments, that $\Phi$ 
is a non-vanishing function because, thanks to property~\textbf{H1}, $\Phi(\mu)<0$ at any admissible parameter $\mu$. On the other hand, from \refc{ex2eq6},
\begin{equation}\label{ex2eq5}
\log\left(\frac{L_{31}}{L_{11}}\right)\!(1)=-\frac{(1+c+cq)}{2c}\log c\text{ and }\log\left(\frac{L_{22}}{L_{32}}\right)\!(1)=0.
\end{equation}
Accordingly the substitution of \refc{ex2eq4} and \refc{ex2eq5} in \refc{ex2d2} yields
\begin{align*}
d_2(\mu)&=-\frac{m}{2q}\log\left(\frac{p}{q}\right)+p(m\bb+2np)\Phi(\mu)
-\frac{p(1+c+cq)}{2qc}\log c \\[7pt]
&=\frac{pq+p+q}{2q}\log\left(\frac{p}{q}\right)+p\big(2pqa-(pq-p+q)\bb\big)\Phi(\mu)
-\frac{p(1+c+cq)}{2qc}\log c,
\end{align*}
where in the first equality we set the values of the hyperbolicity ratios given in \refc{ex2eq0} and in the second one the expressions of $m$ and $n$ defined in \refc{ex2eq3}. Observe at this point, see \refc{ex2eq1}, that $d_1(\mu)=0$ if, and only if, $p=cq.$ 
Moreover the two logarithmic summands in the above expression of $d_2(\mu)$ cancel each other after the substitution $p=cq$, so that
\[
 d_2(\mu)\big|_{p=cq}=cq^2\big(2cqa-(cq-c+1)\bb\big)\Phi(a,\bb,c,cq,q).
\]
Thus, by the Weierstrass Division Theorem (see \cite{Greuel,Krantz}), there exists an analytic function~$\kappa_1$ such that 
\begin{equation}\label{d2div}
d_2(\mu)=d_1(\mu)\kappa_1(\mu)+\big(2cqa-(cq-c+1)\bb\big)\kappa_2(\mu),
\end{equation}
where $\kappa_2(\mu)\!:=cq^2\Phi(\mu)$ is a unity in the admissible set. This expression shows that if we take an admissible parameter $\mu_0=(a_0,\bb_0,c_0,p_0,q_0)$ such that $p_0-c_0q_0=0$ and $2c_0q_0a_0-(c_0q_0-c_0+1)\bb_0\neq 0$ then $d_2(\mu_0)\neq 0$, which by $(c)$ in \teoc{thmA} implies that $\mathrm{Cycl}\big((\Gamma,X_{\mu_0}),X_\mu\big)\leqslant 1.$ The fact that $\mathrm{Cycl}\big((\Gamma,X_{\mu_0}),X_\mu\big)=1$ follows by applying $(b)$ in \teoc{thmA} because $\nabla d_1(\mu_0)$ is not the zero vector, see \refc{ex2eq1}. This proves the validity of the assertion in $(b)$.

To show $(c)$ we take any $\mu_0=(a_0,\bb_0,c_0,p_0,q_0)$ satisfying $p_0-c_0q_0=2c_0q_0a_0-(c_0q_0-c_0+1)\bb_0=0$. Then one can verify that the function
\[
H(x,y)=\frac{q_0(x^2+c_0(y^2+1))-\bb_0 xy}{(xy^{c_0})^\frac{2}{c_0q_0+c_0+1}}
\]
is a first integral of $X_{\mu_0}$, which is clearly analytic on the whole first quadrant. (For reader's convenience let us mention that we found this first integral looking for an integrating factor of the form $x^ry^s$ with $r,s\in\R.$) 
Thus, since the determinant of the Jacobian of $X_{\mu_0}$ at $Q_{\mu_0}$ is strictly positive, we can assert that it is a center. A straightforward application of the Poincaré-Bendixson theorem shows that $\Gamma$ is the outer boundary of its period annulus, which fills the first quadrant. This proves that
if $p_0-c_0q_0=2c_0q_0a_0-(c_0q_0-c_0+1)\bb_0=0$ then the displacement map $\mathscr D(\,\cdot\,;\mu_0)$ of $X_{\mu_0}$ along~$\Gamma$ is identically zero. The converse follows by \teoc{prop1} noting that $d_1(\mu)=d_2(\mu)=0$ if, and only if, $p-cq=2cqa-(cq-c+1)\bb=0$. 

It only remains to be proved that $\mathrm{Cycl}\big((\Gamma,X_{\mu_0}),X_\mu\big)=1$. This follows verbatim the proof of the same fact in assertion $(c)$ of \teoc{ex1} and so we shall omit the details for the sake of shortness. Indeed, by~$(b)$ in \teoc{prop1} we have that if $\mu\approx\mu_0$ then
\[
 \mathscr D(s;\mu)=(b_1(\mu)\omega(s;\alpha(\mu))+b_2(\mu)+r(s;\mu))s^{c},
\]
where $\alpha=c-p/q$ and $r\in\F_\ell^\infty(\mu_0)$ with $\ell\in(0,\min(1,p))$. Here we use that the three hyperbolicity ratios are $\lambda_1=1/q$, $\lambda_2=p$ and $\lambda_3=1/c$. We also have that $(b_1)=(d_1)$ and 
$(b_1,b_2)=(d_1,d_2)$ over the local ring $\R\{\mu\}_{\mu_0}$. Therefore, if $\mu=(a,\bb,c,p,q)$ verifies $p-cq=2cqa-(cq-c+1)\bb=0$ then $b_1(\mu)=b_2(\mu)=0$ and $r(s;\mu)\equiv 0$. Moreover, since $\nabla d_1(\mu_0)$ and $\nabla d_2(\mu_0)$ are linearly independent, see~\refc{ex2eq1} and~\refc{d2div}, this is also the case of 
$\nabla b_1(\mu_0)$ and $\nabla b_2(\mu_0)$. We can thus take three linear functions, say $b_3(\mu)$, $b_4(\mu)$ and $b_5(\mu)$, such that $\nu=\Psi(\mu)\!:=\big(b_1(\mu),b_2(\mu),b_3(\mu),b_4(\mu),b_5(\mu)\big)$ is a local analytic change of coordinates in a neighbourhood of $\mu=\mu_0$ with $\Psi(\mu_0)=0_5\!:=(0,0,0,0,0).$ Notice then that 
$\Psi$ maps $\{p-cq=2cqa-(cq-c+1)\bb=0\}$ to $\{\nu_1=\nu_2=0\}$ and, moreover,
\[
 \mathscr R_1(s;\nu):=s^{-c}\mathscr D(s;\mu)\Big|_{\mu=\Psi^{-1}(\nu)}=\nu_1\omega(s;\hat\alpha)+\nu_2+\hat r(s;\nu),
\]
where $\hat\alpha=\hat\alpha(\nu)\!:=\alpha(\Psi^{-1}(\nu))$ and $\hat r(s;\nu)\!:=r(s;\hat \Psi^{-1}(\nu))\in\F_\ell^\infty(0_5)$.
Due to $\hat r(s;0,0,\nu_3,\nu_4,\nu_5)\equiv 0$, by applying \cite[Lemma 4.1]{MV22} we can write the remainder as $\hat r(s;\nu)=\nu_1h_1(s;\nu)+\nu_2h_2(s;\nu)$
with $h_i\in\F_{\ell}^\infty(0_5)$ and, consequently,
\begin{equation*}
 \mathscr R_1(s;\nu)=\nu_1\big(\omega(s;\hat\alpha)+h_1(s;\nu)\big)+\nu_2\big(1+h_2(s;\nu)\big).
\end{equation*}
From this expression we conclude that there exists $s_0>0$ and an open neighbourhood~$U$ of $\nu=0_5$ such that
$\mathscr R_1(s;\nu)$ that has at most one zero on $(0,s_0)$, counted with multiplicities, for all $\nu=(\nu_1,\nu_2,\nu_3,\nu_4,\nu_5)\in U$ with $\nu_1^2+\nu_2^2\neq 0$, which implies that $\mathrm{Cycl}\big((\Gamma,X_{\mu_0}),X_\mu\big)\leqslant 1$. The proof of this follows exactly as we argue to show the same fact in \teoc{ex1}, cf.~\refc{ex1eq10}, and it is omitted for brevity. Finally the fact that $\mathrm{Cycl}\big((\Gamma,X_{\mu_0}),X_\mu\big)\geqslant 1$ follows taking $\mu_1\approx\mu_0$ with $p_1-c_1q_1=0$ and $2c_1q_1a_1-(c_1q_1-c_1+1)\bb_1\neq 0$, and applying the assertion in $(b)$. This completes the proof of the result.
\end{prooftext}

\begin{obs}
In order to prove \teoc{ex2} it is only necessary to compute the functions $d_1$ and $d_2$ in \teoc{thmA}, which give the conditions for cyclicity 0 and 1, respectively. Let us explain that, as a matter of fact, we computed the function $d_3$ as well, realizing that it vanishes when $d_1=d_2=0.$ It was this fact that lead us to investigate if the return map along the polycycle is the identity in that case. For completeness let us explain succinctly the computations that involve the obtention of $d_3$ for the Kolmogorov's family considered in \teoc{ex2}. Recall, see $(e)$ in \teoc{thmA}, that
\[
d_3(\mu)\!:=\hat M_3(\lambda_3,1)L_{11}(1)-\hat M_1\big(1/\lambda_1,1\big)L_{31}(1).
\]
In this case, cf. \refc{ex2eq6}, we have that
\[
 L_{11}(u)=(1+cu^2)^{\frac{1+(q+1)c}{2c}}\text{ and }L_{31}(u)=\left(1+u^2/c\right)^{\frac{1+(q+1)c}{2c}}
\]
and then, from the definition in \refc{def_M},
\begin{align*}
M_1(u)&=(1+cu^2)^{\frac{1+(q-3)c}{2c}}\big(aq+\bb+(c\bb-(c+1)a)u^2\big)\\
\intertext{and}
M_3(u)&=-u\left(1+u^2/c\right)^{\frac{1+(q-3)c}{2c}}\big((aq+\bb)u^2+c\bb-(c+1)a\big)/c^2.
\end{align*}
In order to proceed with the computation of $\hat M_1\big(\frac{1}{\lambda_1},1\big)$ and $\hat M_3\big(\lambda_3,1\big)$  we note that if $J(x;\eta,r)\!:=(1+\eta x^2)^r$ with $\eta>0$ and $r\in\R$ then 
\[
 \hat J(\alpha,1;\eta,r)=\int_0^1(1+\eta x^2)^rx^{-\alpha-1}dx=-\frac{1}{\alpha}\,{}_2F_1(-r,-\alpha/2;1-\alpha/2;-\eta)\text{ for all $\alpha<0$},
\]
where in the first equality we apply $(b)$ in \propc{L8} with $k=0$ and in the second one we use the equality in \cite[15.3.1]{AbramovitzStegun} to express the definite integral as a hypergeometric function. In principle the above equality is only true provided that $\alpha<0$. However its validity can be extended to any $\alpha\notin\N$ thanks to the meromorphic properties of the functions ${}_2F_1$ and $\hat J$ stablished, respectively, by \cite[Lemma B.2]{MV22} and~$(d)$ in \propc{L8}. Consequently, thanks to this observation and applying twice the above formula, we get
\[
\hat M_1(1/\lambda_1,1)=-\frac{aq+\bb}{q}\,\varphi_1(c,q)+\frac{c\bb-(c+1)a}{2-q}\,\varphi_2(c,q),
\]
where
\begin{align*}
\varphi_1(c,q)&\!:={}_2F_1\big((3-q)/2-1/(2c),-q/2;1-q/2;-c\big),\\[5pt]
\varphi_2(c,q)&\!:={}_2F_1\big((3-q)/2-1/(2c),1-q/2;2-q/2;-c\big).
\end{align*}
Here we also use that if $h=f+g$ then $\hat h_\alpha=\hat f_\alpha+\hat g_\alpha$ and that if $f(x)=x^ng(x)$ then $\hat f_\alpha(x)=x^n\hat g_{\alpha-n}(x),$ see \cite[Corollary B3]{MV21}. Similarly
\[
\hat M_3(\lambda_3,1)=\frac{aq+\bb}{c(1-3c)}\,\varphi_3(c,q)+\frac{c\bb-(c+1)a}{c(1-c)}\,\varphi_4(c,q),
\]
where
\begin{align*}
\varphi_3(c,q)&\!:={}_2F_1\big((3-q)/2-1/(2c),3/2-1/(2c);5/2-1/(2c);-1/c\big),\\[5pt]
\varphi_4(c,q)&\!:= {}_2F_1\big((3-q)/2-1/(2c),1/2-1/(2c);3/2-1/(2c);-1/c\big).
\end{align*}
In the proof of \teoc{ex2} we show that $d_1(\mu)=d_2(\mu)=0$ if, and only if, $\mu=(a,b,c,p,q)$ verifies 
$p=cq$ and $a=\frac{\bb(1-c+cq)}{2cq}$. Long but easy computations show that, under these two conditions, $d_3(\mu)=0$ if, and only if,
\[
\frac{q}{1-3c}\varphi_3(c,q)-\varphi_4(c,q)+c^{-\frac{1+c+cq}{2c}}\left(\varphi_1(c,q)+\frac{c-1}{q-2}\varphi_2(c,q)\right)=0.
\]
This is an equation for $b,$ $c$ and $q$ that involves four hypergeometric functions. Surprisingly enough it turns out, by applying the formula in \cite[15.3.7]{AbramovitzStegun}, that the function on the left hand side of the above equation is identically zero. In other words, $d_1(\mu_0)=d_2(\mu_0)=0$ implies $d_3(\mu_0)=0.$
\end{obs}

\begin{proclama}{Acknowledgements.}
The authors want to thank Joan Torregrosa, who kindly lend us his Maple~\cite{Maple} procedure to compute the focus quantities that we use to prove $(e)$ in \teoc{ex1}.
\end{proclama}

\bibliographystyle{plain}

\end{document}